\documentclass[ejs]{imsart}

\newcommand{\Var}{\mbox{Var}}

\newcommand{\ind}{\stackrel{\mathrm{ind}}{\sim}}

\usepackage{amsthm}
\usepackage{amsmath}
\usepackage{amssymb, amsfonts}
\usepackage{natbib}
\usepackage{bm}
\usepackage{graphicx}
\usepackage{color}
\usepackage{floatrow}
\usepackage{units}
\usepackage{url}
\usepackage{verbatim}
\usepackage{geometry}
\usepackage{slashbox}
\usepackage{epstopdf}
\usepackage{booktabs,caption,fixltx2e}
\usepackage[flushleft]{threeparttable}
\usepackage{rotating}
\usepackage{multirow}
\usepackage{physics}
\usepackage{hyperref}
\usepackage{nameref}
\usepackage{mathptmx} 
\usepackage{placeins} 

\newcommand{\mbf}[1]{\mathbf{#1}}

\newtheorem{theorem}{Theorem}[section]
\newtheorem{lemma}[theorem]{Lemma}

\makeatletter
\let\orgdescriptionlabel\descriptionlabel
\renewcommand*{\descriptionlabel}[1]{%
  \let\orglabel\label
  \let\label\@gobble
  \phantomsection
  \edef\@currentlabel{#1}%
  \let\label\orglabel
  \orgdescriptionlabel{#1}%
}
\makeatother

\arxiv{math.PR/0000000}
\begin{document}

\begin{frontmatter}
\title{Bayesian Estimation Under Informative Sampling}
\runtitle{Bayesian Estimation Under Informative Sampling}
\thankstext{T1}{U.S. Bureau of Labor Statistics, 2 Massachusetts Ave. N.E, Washington, D.C. 20212 USA}

\begin{aug}
\author{\fnms{Terrance D. } \snm{Savitsky}\ead[label=e1]{Savitsky.Terrance@bls.gov}}
\and
\author{\fnms{Daniell} \snm{Toth}\ead[label=e2]{Toth.Daniell@bls.gov}}

\runauthor{T. D. Savitsky et al.}

\affiliation{U.S. Bureau of Labor Statistics}
\address{2 Massachusetts Ave. N.E, Washington, D.C. 20212 USA\\
\printead{e1}\\
\printead{e2}}
\end{aug}

\maketitle
\begin{abstract}
Bayesian analysis is increasingly popular for use in social science and other application areas where the
data are observations from an informative sample.
An informative sampling design leads to inclusion probabilities that are correlated with the response variable of interest.  Model inference performed on the observed sample taken from the population will be biased for the population generative model under informative sampling since the balance of information in the sample data is different from that for the population.  Typical approaches to account for an informative sampling design under Bayesian estimation are often difficult to implement because they require re-parameterization of the hypothesized generating model, or focus on design, rather than model-based, inference. We propose to construct a pseudo-posterior distribution that utilizes sampling weights based on the marginal inclusion probabilities to exponentiate the likelihood contribution of each sampled unit, which weights the information in the sample back to the population.  Our approach provides a nearly automated estimation procedure applicable to any model specified by the data analyst for the population and retains the population model parameterization and posterior sampling geometry. We construct conditions on known marginal and pairwise inclusion probabilities that define a class of sampling designs where $L_{1}$ consistency of the pseudo posterior is guaranteed.   We demonstrate our method on an application concerning the Bureau of Labor Statistics Job Openings and Labor Turnover Survey.
\end{abstract}


\end{frontmatter}

\section{Introduction} \label{motivation}

Bayesian formulations are increasingly popular for modeling hypothesized distributions
with complicated dependence structures.
Their popularity stems from the ease of capturing this dependence by employing models with random effects parameters with a hierarchical construction that regulates the borrowing of information for estimation.
Latent parameters are often used in the model to permit flexibility in the estimation of the dependencies among the observations \citep{dunson2010nonparametric}.  In social science applications, utilization of latent parameters may be useful for making inference about intrinsic belief states of people from their observed actions(see for example, \citet{RSSC:RSSC12049})  Other application areas in which latent parameters may be employed include, engineering and natural science, which use them to parameterize elements of
an evolving process.

Data used in these type of applications are often acquired through a complex sample design,
resulting in probabilities of inclusion that are associated with the 
variable of interest.
This association could result in an observed data set consisting of units that are not independent and identically distributed.  A sampling design that produces a correlation between selection probabilities and observed values is referred to as informative.  Failure to account for this dependence caused by the sampling design could bias estimation of parameters that index the joint distribution hypothesized to have generated the population (\citep{holt:1980}).

\subsection{Examples}\label{examples}
We next outline some examples of survey instruments that employ informative sampling designs and associated inferential goals for models estimated on observed samples realized from these surveys.

\textit{Example 1:} The Survey of Occupational Illnesses and Injuries (SOII) is administered to U.S. business establishments by the U.S. Bureau of Labor Statistics (BLS), in partnership with individual states, in order to capture workplace induced injuries and illnesses.  A stratified sampling design is used where strata are indexed by state-industry-size-injury rate.  Strata containing establishments that historically express higher injury rates are assigned higher sample inclusion probabilities.  The resulting sample will contain a larger proportion of establishments that express higher injury rates than the population, as a whole.  States desire to perform regression modeling with variable selection to discover the root causes that predict illnesses and injuries among the population of establishments, estimated from the observed sample.  The model-estimated coefficients from the sample will be biased absent correction for over-representation of establishments that tend to express relatively high injury rates.

\textit{Example 2:} The Current Establishment Statistics (CES) is a BLS survey of U.S. business establishments that collects employment count data across states and industries under a stratified sampling design with strata indexed by the number of employees in each establishment.  Strata containing relatively larger establishments are assigned higher inclusion probabilities than those which hold establishments with relatively fewer employees.  The distribution of employment totals in the observed sample of establishments will be skewed towards relatively larger values as compared to the population of establishments. An important area of modeling inference is to understand industry-indexed differences in monthly employment trends and correlations among industries in the population.   We would use a mixed effects model, parameterized with random effects indexed by industry and month.   Estimation of the population distribution under our model from the observed sample will be biased absent some correction for the skewness in the sample towards larger-sized establishments.

\textit{Example 3:} BLS collects establishment-indexed employment totals in \emph{both} the Quarterly Census of Employment and Wages (QCEW) \emph{and} the CES survey.  CES survey participants also provide submissions to the QCEW, such that their reported monthly employment totals for an overlapping time period of interest should be equal between the two instruments, but they are not for approximately $10000$ establishments, indicating one or more employment count submission errors for those respondents.  A response variable of interest, termed the ``error time series", was created by taking the absolute value of the difference in reported employment totals among the $10000$ establishments for each month over a $12$ month period.  A ``response analysis survey" (RAS) of approximately $2000$ establishments was taken from this population with the goal to understand the process drivers for committing errors so that BLS may target resources to establishments that mitigate them.  The modeling focus is to identify probabilistic clusters of establishments with similar error patterns over the $12$ month period and to examine the process by which establishments in each cluster construct their data submissions to BLS.  The RAS survey design stratified the population of $10000$ establishments based on phenomena of interest expressed in portions of each time series; for example, a big jump in the reported difference at year-end may indicate establishments who count checks that include regular pay and bonuses for each employee, instead of counting employees.  Higher inclusion probabilities were assigned to those strata expressing phenomena of relatively greater interest to BLS researchers.  Modeling the number of and memberships in probabilistic clusters of error patterns expressed in the population from the RAS sample may be biased because the proportions of error patterns expressed in the sample are designed to be different from the population.

\textit{Example 4:} The Current Expenditure (CE) survey is administered to U.S. households by BLS for the purpose of determining the amount of spending for a broad collection of goods and service categories and it serves as the main source used to construct the basket of goods later used to formulate the Consumer Price Index.  The CE employs a multi-stage sampling design that draws clusters of core-based statistical areas (CBSAs), such as metropolitan and micropolitan areas, from which Census blocks and, ultimately, households are sampled.  Economists desire to model the propensity or probability of purchase for a variety of goods and services.  The balance of sampled clusters may not be reflective of those in the population; for example, if particularly high income ares are included in the sample.  So inference on purchase propensities for the population made from the observed sample will be biased absent correction for the informative sampling design.

\textit{Example 5:} BLS administers the Job Openings and Labor Turnover survey (JOLTS) to business establishments with the focus to measure labor market dynamics by reporting the number of job openings, hires and separations, which is a leading indicator for employment trends.  The sampling design assigns larger inclusion probabilities to establishments with relatively more employees because larger establishments drive the variance in the reported statistics.  Our modeling goals are to understand differences in labor force dynamics based on employment ownership (e.g., private, public) and region as part of imputing missing values with respect to the population generating distribution.  As with the CES sampling design, however, our sample will tend to over-represent relatively larger-sized establishments, so that inference and imputation using the sample will be biased for the population.  We develop a multivariate count data population generating model in Section~\ref{application}, where we illustrate the resulting estimation bias from failure to account for the correlations between assigned inclusion probabilities and the response variables of interest for our sample.

The target audience for this article are data analysts who wish to perform some distributional inference
using data obtained from an informative sample design on a population using a model they specify, $p\left(y_{i}\vert \bm{\lambda}\right),~\bm{\lambda} \in \Lambda$, for density, $p$.
We discuss, in the next section, how the limited literature on this topic does not adequately provide a
general method for making distributional inference on a population while
adjusting for the unequal probabilities of selection.

In this article, we propose an approach that replaces the likelihood with the ``pseudo" likelihood \citep{chambers2003analysis}, $p\left(y_{i}\vert \delta_{i} = 1,\bm{\lambda}\right)^{w_{i}}$, using sampling weight, $w_{i} \propto 1/\pi_{i}$.
This re-weights the likelihood contribution for each observed unit with intent to re-balance the information in the observed sample to approximate the balance of information in the target finite population;
correcting for the informativeness.
We show that the proposed method for Bayesian estimation on complex sample data allows for asymptotically consistent inference on any population-generating model specified by the data analyst.

Additionally, this method does not require information about the complex design, other than the probabilities of selection, or about the full population, other than the observed data.
We believe this makes the method applicable to more situations.
Indeed, it is often the case that the data analyst does not have access to the full design information or auxiliary variables on the population,
$z_{1},\ldots,z_{N}$, used to assign the probabilities of selection $\pi_{1},\ldots,\pi_{N}.$
However, it is common for the probabilities of selection for the units in the sample, $\pi_{1},\ldots,\pi_{n},$
to be provided with the observed sample data.


\subsection{Review of  Methods to Account for Informative Sampling}
One current approach is to account for the informativeness by parameterizing the sampling design into the model \citep{Litt:mode:2004}.
Parameterizing even a simple informative design is often difficult to accomplish and may disrupt desired inference by requiring a change to the underlying population model parameterization.  The analyst in Example $3$, above, desires to perform inference on an \emph{a priori} unknown clustering of sampled units with their population model for data acquired under a stratified sampling design.  Specifying random effects to be indexed by strata will likely conflict with the identification and composition of inferred clusters. Further, the data analyst may not have access to the sampling design, but only indirect information in form of sampling weights.  Lastly, the analyst is sometimes required to impute the unobserved units in the finite population, which may be computationally infeasible.

Another approach incorporates the sampling weights into inference about the population, as is our intent, but requires a particular form for the likelihood that does not allow the analyst to impose their \emph{own} population model formulation of inferential interest. For example, \citet{dong:2014} specifies an empirical likelihood, while \citet{kunihama:2014} constructs a non-parametric mixture for the likelihood and \citet{wu:2010} uses a sampling-weighted (pseudo) empirical likelihood.  All of these approaches impose Dirichlet distribution priors for the mixture components with hyperparameters specified as a function of the first-order sampling weights. \citet{si2015} regress the response variable on a Gaussian process function of the weights for sampling designs where sub-groups of sampled units have equal weights (e.g., a stratified sampling design).   These approaches are designed for inference about simple mean and total statistics, rather than inference for parameters that characterize an analyst-specified population model that is the focus for our proposed method.

One method that uses a plug-in estimator, as do we in our method, is to construct a joint likelihood of the population distribution and sample inclusion in a simple logistic regression model \citep{Male:Davi:Cao:mode:1999}.
This allows one to analytically marginalize over the parameters indexed by the non-sampled units.
This approach is limited in application to a class of simple population models that permit analytic integration and may not be applied to more general classes of Bayesian models for the population that we envision in development of our approach.

Perhaps the most general Bayesian approach constructs models to co-estimate parameters for conditional expectations of inclusion probabilities jointly with the population-generating model parameters at each level of a hierarchical construction \citep{Pfef:DaS:DoN:mult:2006}.
This formulation is fully Bayesian such that it accounts for all sources of uncertainty in population generation and inclusion of units,
but requires a custom implementation of an MCMC sampler for each specified population model, such as their simple two-level linear regression model.   The implementations may increase the complexity of the specified model and reduce the quality of posterior mixing in the MCMC, so that they are suitable for relatively simple population probability models.

The method we propose is intended to allow Bayesian inference from \emph{any} population model that may be specified by the the data analyst under an informative sampling design, unlike the alternative methods. It provides asymptotically unbiased estimation using \emph{only} the distribution for the observed sample units and normalized H\'{a}jek-like sampling weights.  The ``plug-in" type method accounts for the informative sampling design by raising the likelihood contribution of each sampled observation to the power of their associated sampling weight.  The implementation of the plug-in procedure for Bayesian estimation multiplies the sampling weight into each full conditional log-posterior density. This can then sampled in the typical sequential scan MCMC.

Unlike these other methods that are prominent in the literature, this method:
1.) does not impose a population model (implicitly or explicitly), unlike the most recently-developed methods \citep{dong:2014,kunihama:2014,wu:2010,si2015}; 2.) requires only the sampling weights and does not require parameterizing the sampling design unlike \citet{Litt:mode:2004}; 3.) does not require a customized MCMC sampling procedure unlike \citet{Pfef:DaS:DoN:mult:2006}, so can be done automatically; 4). does not require imputing the non-sampled units in the finite population.  Our data application and estimation model in the sequel are intended to be representative of common problems for Bayesian inference, and the application data are not readily estimated with these other methods that account for informative sampling.


We formulate the pseudo-posterior density as sampling weight-adjusted plug-in from which we conduct model inference about the population under a dependent, informative sampling design in Section~\ref{methods}.
Conditions are constructed that guarantee a frequentist $L_{1}$ contraction of the pseudo posterior distribution on the true generating distribution in Section~\ref{theory}.  We make an application of the pseudo posterior estimator to construct a regression model for count data using a dataset of monthly job hires and separations collected by the U.S. Bureau of Labor Statistics in Section~\ref{application}.  We reveal large differences for parameter estimates between incorporation versus ignoring the sampling weights.  This section also includes a simulation study that compares the pseudo posterior estimated on the observed sample to the posterior estimated on the entire finite population.  The paper concludes with a discussion in Section~\ref{discussion}.  The proofs for the main result, along with two enabling results are contained in an Appendix.

\section{Method to account for Informative Sampling}\label{methods}
We begin by constructing the pseudo likelihood and associated pseudo posterior density under any analyst-specified prior formulation on the model, $\bm{\lambda} \in \Lambda$.

\subsection{Pseudo Posterior} \label{bonnery}
Suppose there exists a Lebesgue measurable population-generating density, $\pi\left(y\vert\bm{\lambda}\right)$, indexed by parameters, $\bm{\lambda} \in \Lambda$. Let $\delta_{i} \in \{0,1\}$ denote the sample inclusion indicator for units $i = 1,\ldots,N$ from the population under sampling without replacement.  The density for the observed sample is denoted by, $\pi\left(y_{o}\vert\bm{\lambda}\right) = \pi\left(y\vert \delta_{i} = 1,\bm{\lambda}\right)$, where ``$o$" indicates ``observed".

The plug-in estimator for posterior density under the analyst-specified model for $\bm{\lambda} \in \Lambda$ is
\begin{equation}
\hat{\pi}\left(\bm{\lambda}\vert \mathbf{y}_{o},\tilde{\mathbf{w}}\right) \propto \left[\mathop{\prod}_{i = 1}^{n}p\left(y_{o,i}\vert \bm{\lambda}\right)^{\tilde{w}_{i}}\right]\pi\left(\bm{\lambda}\right), \label{pseudolike}
\end{equation}
where
$\mathop{\prod}_{i=1}^{n}p\left(y_{o,i}\vert \bm{\lambda}\right)^{\tilde{w}_{i}}$ denotes the pseudo likelihood for observed sample responses, $\mathbf{y}_{o}$. The joint prior density on model space assigned by the analyst is denoted by $\pi\left(\bm{\lambda}\right)$.  This pseudo likelihood employs sampling weights, $\{\tilde{w}_{i} \propto 1/\pi_{i}\}$, constructed to be inversely proportional to unit inclusion probabilities.  Each sampling weight assigns the relative importance of the likelihood contribution for each sample observation to approximate the likelihood for the population.  We use $\hat{\pi}$ to denote the noisy approximation to posterior distribution, $\pi$, and we make note that the approximation is based on the data, $\mathbf{y}_{o}$ , and sampling weights, $\{\tilde{\mathbf{w}}\}$, confined to those units \emph{included} in the sample, $S$.

The total estimated posterior variance is regulated by the sum of the sampling weights.  We define unnormalized weights, $\{w_{i} = 1/\pi_{i}\}$, and subsequently normalize them, $\tilde{w}_{i} = \frac{w_{i}}{\frac{\sum w_{i}}{n}},~i = 1,\ldots,n$, to sum to the sample size, $n$, the asymptotic units of information in the sample.  Incorporation of the sampling weights to formulate the pseudo posterior estimator is expected to increase the estimated parameter posterior variances relative to the (unweighted) posterior estimated on a simple random (non-informative) sample because the weights encode the uncertainty with which samples represent the finite population under repeated sampling. This increase in estimated posterior variance may be partly or wholly offset to the extent that the informative sampling design is more efficient than simple random sampling; for example, a stratified sampling design that takes simple random samples within each stratum may produce samples that provide better coverage of the population.  Although our method utilizes the weights as a ``plug-in", rather than imposing a prior, \citet{Pfef:Sver:infe:2009} use Bayes rule to demonstrate one may replace the weights with their conditional expectation given the observed response to correct for informative sampling.  Replacing the raw weights with their conditional expectation given the observed response may serve to reduce the total variation attributed to weighting (and the resulting posterior uncertainty) in the case where the actual sampled observations express information in different proportions than intended in the sampling design.  Even though the conditional \emph{distribution} of the weights given the response is generally different for the observed sample than for the population, nevertheless their conditional \emph{expectations} are equal.

\section{Pseudo Posterior Consistency}\label{theory}
We formulate a pseudo posterior distribution in this section and specify conditions under which it contracts on the true generating distribution in $L_{1}$.  Let $\nu \in \mathbb{Z}^{+}$ index a sequence of finite populations, $\{U_{\nu}\}_{\nu=1,\ldots,N_{\nu}}$, each of size, $\vert U_{\nu} \vert = N_{\nu}$, such that $N_{\nu} < N_{\nu^{'}},\text{ for}~\nu < \nu^{'}$, so that the finite population size grows as $\nu$ increases.  Suppose that $\mbf{X}_{\nu,1},\ldots, \mbf{X}_{\nu,N_{\nu}}$ are independently distributed according to some unknown distribution $P,$ (with density, $p$) defined on the sample space, $\left(\mathcal{X},\mathcal{A}\right).$  If $\Pi$ is a prior distribution on the model space, $\left(\mathcal{P}, \mathcal{C}\right)$ to which $P$ is known to belong, then the posterior distribution is given by
\begin{equation}
\label{postpop}
\Pi\left(B\vert \mbf{X}_{1},\ldots,\mbf{X}_{N_{\nu}}\right) = \frac{\mathop{\int}_{P \in B}\mathop{\prod}_{i=1}^{N_{\nu}}\frac{p}{p_{0}}(\mbf{X}_{i})d\Pi(P)}{\mathop{\int}_{P \in \mathcal{P}}\mathop{\prod}_{i=1}^{N_{\nu}}\frac{p}{p_{0}}(\mbf{X}_{i})d\Pi(P)},
\end{equation}
for any $B \in \mathcal{C}$, where we refer to $\{\mbf{X}_{\nu,i}\}_{i=1,\ldots,N_{\nu}}$ as $\{\mbf{X}_{i}\}_{i=1,\ldots,N_{\nu}}$ for readability when the context is clear.

\citet{ghosal2007} study the rate at which this posterior distribution converges to the assumed true (and fixed) generating distribution $P_{0}.$  They prove, under certain conditions on the model space, $\mathcal{P}$, and the prior distribution, $\Pi$, 
that in $P_{0}-$probability, the posterior distribution concentrates on an arbitrarily small neighborhood of
$P_0$ as $N_{\nu} \uparrow \infty.$

The observed data on which we focus is not the entire finite population, $\mbf{X}_{1},\ldots, \mbf{X}_{N_{\nu}},$ but rather a sample, $\mbf{X}_{1},\ldots, \mbf{X}_{n_{\nu}},$ with $n_{\nu} \leq N_{\nu}$, drawn under a sampling design distribution applied to the finite population under which each unit, $i \in (1,\ldots,N_{\nu})$, is assigned a probability of inclusion in the sample.  These unit inclusion probabilities are constructed to depend on the realized finite population values, $\mbf{X}_{1},\ldots, \mbf{X}_{N_{\nu}}$, at each $\nu$.

\subsection{Pseudo Posterior Distribution}\label{pseudop}
A sampling design is defined by placing a \emph{known} distribution on a vector of inclusion indicators, $\bm{\delta}_{\nu} = \left(\delta_{\nu 1},\ldots,\delta_{\nu N_{\nu}}\right)$, linked to the units comprising the population, $U_{\nu}$.  The sampling distribution is subsequently used to take an \emph{observed} random sample of size $n_{\nu} \leq N_{\nu}$.
Our conditions needed for the main result employ known marginal unit inclusion probabilities, $\pi_{\nu i} = \mbox{Pr}\{\delta_{\nu i} = 1\}$ for all $i \in U_{\nu}$ and the second-order pairwise probabilities, $\pi_{\nu ij} = \mbox{Pr}\{\delta_{\nu i} = 1 \cap \delta_{\nu j} = 1\}$ for $i,j \in U_{\nu}$, which are obtained from the joint distribution over $\left(\delta_{\nu 1},\ldots,\delta_{\nu N_{\nu}}\right)$. The dependence among unit inclusions in the sample contrasts with the usual $iid$ draws from $P$.  We denote the sampling distribution by $P_{\nu}$.

Under informative sampling, the marginal inclusion probabilities, $\pi_{\nu i} = P\{\delta_{\nu i} =1\},~ i \in\left(1,\ldots,N_{\nu}\right)$, are formulated to depend on the finite population data values, $\mathbf{X}_{N_{\nu}} = \left(\mbf{X}_{1},\ldots,\mbf{X}_{N_{\nu}}\right)$.
Since the resulting balance of information would be different in the sample, the posterior distribution for $\left(\mbf{X}_{1}\delta_{\nu 1},\ldots,\mbf{X}_{N_{\nu}}\delta_{\nu N_{\nu}}\right)$, that we employ for inference about $P_{0}$, is not equal to that of Equation~\ref{postpop}.

Our task is to perform inference about the population generating distribution, $P_{0}$, using the observed data taken under an informative sampling design.  We account for informative sampling by ``undoing" the sampling design with the weighted estimator,
\begin{equation}
p^{\pi}\left(\mbf{X}_{i}\delta_{\nu i}\right) := p\left(\mbf{X}_{i}\right)^{\delta_{\nu i}/\pi_{\nu i}},~i \in U_{\nu},
\end{equation}
which weights each density contribution, $p(\mbf{X}_{i})$, by the inverse of its marginal inclusion probability. This construction re-weights the likelihood contributions defined on those units randomly-selected for inclusion in the observed sample ($\{i \in U_{\nu}:\delta_{\nu i} = 1\}$) to \emph{approximate} the balance of information in $U_{\nu}$.  This approximation for the population likelihood produces the associated pseudo posterior,
\begin{equation}\label{inform_post}
\Pi^{\pi}\left(B\vert \mbf{X}_{1}\delta_{\nu 1},\ldots,\mbf{X}_{N_{\nu}}\delta_{\nu N_{\nu}}\right) = \frac{\mathop{\int}_{P \in B}\mathop{\prod}_{i=1}^{N_{\nu}}\frac{p^{\pi}}{p_{0}^{\pi}}(\mbf{X}_{i}\delta_{\nu i})d\Pi(P)}{\mathop{\int}_{P \in \mathcal{P}}\mathop{\prod}_{i=1}^{N_{\nu}}\frac{p^{\pi}}{p_{0}^{\pi}}(\mbf{X}_{i}\delta_{\nu i})d\Pi(P)},
\end{equation}
that we use to achieve our required conditions for the rate of contraction of the pseudo posterior distribution on $P_{0}$.  We recall that both $P$ and $\bm{\delta}_{\nu}$ are random variables defined on the space of measures and possible samples, respectively.  Additional conditions are later formulated for the distribution over samples, $P_{\nu}$, drawn under the known sampling design, to achieve contraction of the pseudo posterior on $P_{0}$.  We assume measurability for the sets on which we compute prior, posterior and pseudo posterior probabilities on the joint product space, $\mathcal{X}\times\mathcal{P}$.  For brevity, we use the superscript, $\pi$, to denote the dependence on the known sampling probabilities, $\{\pi_{\nu i}\}_{i=1,\ldots,N_{\nu}}$; for example, $\displaystyle\Pi^{\pi}\left(B\vert \mbf{X}_{1}\delta_{\nu 1},\ldots,\mbf{X}_{N_{\nu}}\delta_{\nu N_{\nu}}\right) := \Pi\left(B\vert \left(\mbf{X}_{1}\delta_{\nu 1},\ldots,\mbf{X}_{N_{\nu}}\delta_{\nu N_{\nu}}\right),\left(\pi_{\nu 1},\ldots,\pi_{\nu N_{\nu}}\right)\right)$.

Our main result is achieved in the limit as $\nu\uparrow\infty$, under the countable set of successively larger-sized populations, $\{U_{\nu}\}_{\nu \in \mathbb{Z}^{+}}$.  We define the associated rate of convergence notation, $\order{b_{\nu}}$, to denote $\mathop{\lim}_{\nu\uparrow\infty}\frac{\order{b_{\nu}}}{b_{\nu}} = 0$.

\subsection{Empirical process functionals}\label{empirical}
We employ the empirical distribution approximation for the joint distribution over population generation and the draw of an informative sample that produces our observed data to formulate our results.  Our empirical distribution construction follows \citet{breslow:2007} and incorporates inverse inclusion probability weights, $\{1/\pi_{\nu i}\}_{i=1,\ldots,N_{\nu}}$, to account for the informative sampling design,
\begin{equation}
\mathbb{P}^{\pi}_{N_{\nu}} = \frac{1}{N_{v}}\mathop{\sum}_{i=1}^{N_{\nu}}\frac{\delta_{\nu i}}{\pi_{\nu i}}\delta\left(\mbf{X}_{i}\right),
\end{equation}
where $\delta\left(\mbf{X}_{i}\right)$ denotes the Dirac delta function, with probability mass $1$ on $\mbf{X}_{i}$ and we recall that $N_{\nu} = \vert U_{\nu} \vert$ denotes the size of of the finite population. This construction contrasts with the usual empirical distribution, $\mathbb{P}_{N_{\nu}} = \frac{1}{N_{v}}\mathop{\sum}_{i=1}^{N_{\nu}}\delta\left(\mbf{X}_{i}\right)$, used to approximate $P \in \mathcal{P}$, the distribution hypothesized to generate the finite population, $U_{\nu}$.

We follow the notational convention of \citet{Ghosal00convergencerates} and define the associated expectation functionals with respect to these empirical distributions by $\mathbb{P}^{\pi}_{N_{\nu}}f = \frac{1}{N_{\nu}}\mathop{\sum}_{i=1}^{N_{\nu}}\frac{\delta_{\nu i}}{\pi_{\nu i}}f\left(\mbf{X}_{i}\right)$.  Similarly, $\mathbb{P}_{N_{\nu}}f = \frac{1}{N_{\nu}}\mathop{\sum}_{i=1}^{N_{\nu}}f\left(\mbf{X}_{i}\right)$.  Lastly, we use the associated centered empirical processes, $\mathbb{G}^{\pi}_{N_{\nu}} = \sqrt{N_{\nu}}\left(\mathbb{P}^{\pi}_{N_{\nu}}-P_{0}\right)$ and $\mathbb{G}_{N_{\nu}} = \sqrt{N_{\nu}}\left(\mathbb{P}_{N_{\nu}}-P_{0}\right)$.

The sampling-weighted, (average) pseudo Hellinger distance between distributions, $P_{1}, P_{2} \in \mathcal{P}$, $d^{\pi,2}_{N_{\nu}}\left(p_{1},p_{2}\right) = \frac{1}{N_{\nu}}\mathop{\sum}_{i=1}^{N_{\nu}}\frac{\delta_{\nu i}}{\pi_{\nu i}}d^{2}\left(p_{1}(\mathbf{X}_{i}),p_{2}(\mathbf{X}_{i})\right)$, where $d\left(p_{1},p_{2}\right) = \left[\mathop{\int}\left(\sqrt{p_{1}}-\sqrt{p_{2}}\right)^{2}d\mu\right]^{\frac{1}{2}}$ (for dominating measure, $\mu$).
We need this empirical average distance metric because the observed (sample) data drawn from the finite population under $P_{\nu}$ are no longer independent.  The implication is that our result apply to finite populations generated as $inid$ from which informative samples are taken.  The associated non-sampling Hellinger distance is specified with, $d^{2}_{N_{\nu}}\left(p_{1},p_{2}\right) = \frac{1}{N_{\nu}}\mathop{\sum}_{i=1}^{N_{\nu}}d^{2}\left(p_{1}(\mathbf{X}_{i}),p_{2}(\mathbf{X}_{i})\right)$.

\subsection{Main result}\label{results}
We proceed to construct associated conditions and a theorem that contain our main result on the consistency of the pseudo posterior distribution under a class of informative sampling designs at the true generating distribution, $P_{0}$.  Our approach extends the main in-probability convergence result of \citet{ghosal2007} by adding new conditions that restrict the distribution of the informative sampling design.  Suppose we have a  sequence, $\xi_{N_{\nu}} \downarrow 0$ and $N_{\nu}\xi^{2}_{N_{\nu}}\uparrow\infty$  and $n_{\nu}\xi^{2}_{N_{\nu}}\uparrow\infty$ as $\nu\in\mathbb{Z}^{+}~\uparrow\infty$ and any constant, $C >0$,

\begin{description}
\item[(A1)\label{existtests}] (Local entropy condition - Size of model)

        \begin{equation*}
        \mathop{\sup}_{\xi > \xi_{N_{\nu}}}\log N\left(\xi/36,\{P\in\mathcal{P}_{N_{\nu}}: d_{N_{\nu}}\left(P,P_{0}\right) < \xi\},d_{N_{\nu}}\right) \leq N_{\nu} \xi_{N_{\nu}}^{2},
        \end{equation*}
\item[(A2)\label{sizespace}] (Size of space)
        \begin{equation*}
        \displaystyle\Pi\left(\mathcal{P}\backslash\mathcal{P}_{N_{\nu}}\right) \leq \exp\left(-N_{\nu}\xi^{2}_{N_{\nu}}\left(2(1+2C)\right)\right)
        \end{equation*}
\item[(A3)\label{priortruth}] (Prior mass covering the truth)
        \begin{equation*}
        \displaystyle\Pi\left(P: -P_{0}\log\frac{p}{p_{0}}\leq \xi^{2}_{N_{\nu}}\cap P_{0}\left[\log\frac{p}{p_{0}}\right]^{2}\leq \xi^{2}_{N_{\nu}} \right) \geq \exp\left(-N_{\nu}\xi^{2}_{N_{\nu}}C\right)
        \end{equation*}
\item[(A4)\label{bounded}] (Non-zero Inclusion Probabilities)
        \begin{equation*}
        \displaystyle\mathop{\sup}_{\nu}\left[\frac{1}{\displaystyle\mathop{\min}_{i\in U_{\nu}}\pi_{\nu i}}\right] \leq \gamma, \text{  with $P_{0}-$probability $1$.}
        \end{equation*}
\item[(A5)\label{independence}] (Asymptotic Independence Condition)
        \begin{equation*}
        \displaystyle\mathop{\limsup}_{\nu\uparrow\infty} \mathop{\max}_{i \neq j\in U_{\nu}}\left\vert\frac{\pi_{\nu ij}}{\pi_{\nu i}\pi_{\nu j}} - 1\right\vert = \order{N_{\nu}^{-1}}, \text{  with $P_{0}-$probability $1$}
        \end{equation*}
        such that for some constant, $C_{3} > 0$,
        \begin{equation*}
        \displaystyle N_{\nu}\mathop{\sup}_{\nu}\mathop{\max}_{i \neq j \in U_{\nu}}\left[\frac{\pi_{\nu ij}}{\pi_{\nu i}\pi_{\nu j}}\right] \leq C_{3}, \text{  for $N_{\nu}$ sufficiently large.}
        \end{equation*}
\item[(A6)\label{fraction}] (Constant Sampling fraction)
        For some constant, $f \in(0,1)$, that we term the ``sampling fraction",
        \begin{equation*}
        \mathop{\limsup}_{\nu}\displaystyle\biggl\vert\frac{n_{\nu}}{N_{\nu}} - f\biggl\vert = \order{1}, \text{  with $P_{0}-$probability $1$.}
        \end{equation*}
\end{description}
Condition~\nameref{existtests} denotes the logarithm of the covering number, defined as the \emph{minimum} number of balls of radius $\xi/36$ needed to cover $\left\{P\in\mathcal{P}_{N_{\nu}}: d_{N_{\nu}}\left(P,P_{0}\right) < \xi\right\}$ under distance metric, $d_{N_{\nu}}$. This condition restricts the growth in the size of the model space, or as noted by \citet{Ghosal00convergencerates}, the space, $\mathcal{P}_{N_{\nu}}$, must be not too big in order that the condition specifies an optimal convergence rate \citep{wong1995}.  This condition guarantees the existence of test statistics, $\phi_{n_{\nu}}\left(\mbf{X}_{1}\delta_{\nu 1},\ldots,\mbf{X}_{N_{\nu}}\delta_{\nu N_{\nu}}\right) \in\left(0,1\right)$, needed for enabling Lemma~\ref{numerator}, stated in the Appendix, that bounds the expectation of the pseudo posterior mass assigned on the set $\{P\in \mathcal{P}_{N_{\nu}}:d_{n_{\nu}}\left(P,P_{0}\right) \geq \xi_{N_{\nu}}\}$.  Condition~\nameref{priortruth} ensures the prior, $\Pi$, assigns mass to convex balls in the vicinity of $P_{0}$. Conditions~\ref{existtests} and ~\ref{priortruth}, together, define the minimum value of $\xi_{N_{\nu}}$, where if these conditions are satisfied for some $\xi_{N_{\nu}}$, then they are also satisfied for any $\xi > \xi_{N_{\nu}}$. Condition~\nameref{sizespace} allows, but restricts, the prior mass placed on the uncountable portion of the model space, such that we may direct our inference to an approximating sieve, $\mathcal{P}_{N_{\nu}}$.  This sequence of spaces ``trims" away a portion of the space that is not entropy bounded (in condition~\nameref{existtests}).  In practice, trimming the space may usually be performed to ensure the entropy bound.

The next three new conditions impose restrictions on the sampling design and associated known distribution, $P_{\nu}$, used to draw the observed sample data that, together, define a class of allowable sampling designs on which the contraction result for the pseudo posterior is guaranteed.  Condition~\nameref{bounded} requires the sampling design to assign a positive probability for inclusion of every unit in the population because the restriction bounds the sampling inclusion probabilities away from $0$. Since the maximum inclusion probability is $1$, the bound, $\gamma \geq 1$. No portion of the population may be systematically excluded, which would prevent a sample of any size from containing information about the population from which the sample is taken.  Condition~\nameref{independence} restricts the result to sampling designs where the dependence among lowest-level sampled units attenuates to $0$ as $\nu\uparrow\infty$; for example, a two-stage sampling design of clusters within strata would meet this condition if the number of population units nested within each cluster (from which the sample is drawn) increases in the limit of $\nu$.  Such would be the case in a survey of households within each cluster if the cluster domains are geographically defined and would grow in area as $\nu$ increases.  In this case of increasing cluster area, the dependence among the inclusion of any two households in a given cluster would decline as the number of households increases with the size of the area defined for that cluster.  Condition~\nameref{fraction} ensures that the observed sample size, $n_{\nu}$, limits to $\infty$ along with the size of the partially-observed finite population, $N_{\nu}$.

\begin{theorem}
\label{main}
Suppose conditions ~\nameref{existtests}-\nameref{fraction} hold.  Then for sets $\mathcal{P}_{N_{\nu}}\subset\mathcal{P}$, constants, $K >0$, and $M$ sufficiently large,
\begin{align}\label{limit}
&\mathbb{E}_{P_{0},P_{\nu}}\Pi^{\pi}\left(P:d^{\pi}_{N_{\nu}}\left(P,P_{0}\right) \geq M\xi_{N_{\nu}} \vert \mbf{X}_{1}\delta_{\nu 1},\ldots,\mbf{X}_{N_{\nu}}\delta_{\nu N_{\nu}}\right) \leq\nonumber\\
&\frac{16\gamma^{2}\left[\gamma+C_{3}\right]}{\left(Kf + 1 - 2\gamma\right)^{2}N_{\nu}\xi_{N_{\nu}}^{2}} + 5\gamma\exp\left(-\frac{K n_{\nu}\xi_{N_{\nu}}^{2}}{2\gamma}\right),
\end{align}
which tends to $0$ as $\left(n_{\nu}, N_{\nu}\right)\uparrow\infty$.
\end{theorem}

We note that the rate of convergence is injured for a sampling distribution, $P_{\nu}$, that assigns relatively low inclusion probabilities to some units in the finite population such that $\gamma$ will be relatively larger. Samples drawn under a design that expresses a large variability in the sampling weights will express more dispersion in their information similarity to the underlying finite population.  Similarly, the larger the dependence among the finite population unit inclusions induced by $P_{\nu}$, the higher will be $C_{3}$ and the slower will be the rate of contraction.

The separability of the conditions on $\mathcal{P}$ and $\Pi\left(P\right)$, on the one hand, from those on the sampling design distribution, $P_{\nu}$, on the other hand, coupled with the sequential process of taking the observed sample from the finite population reveal that the pseudo posterior, defined on the partially-observed sample from a population, contracts on $P_{0}$ through converging to the posterior distribution defined on each fully-observed population.  We demonstrate this property of the pseudo posterior in a simulation study conducted in Section~\ref{simulation}.  By contrast, if the posterior distribution, defined on each fully-observed finite population, fails to meet conditions~\ref{existtests}, \ref{sizespace} and \ref{priortruth} for the main result from Equation~\ref{limit}, such that it fails to contract on $P_{0}$, then the associated pseudo posterior cannot contract on $P_{0}$, even if the sampling design satisfies conditions~\nameref{bounded}, \nameref{independence} and \nameref{fraction}.

The proof generally follows that of \citet{Ghosal00convergencerates} with substantial modification to account for informative sampling.  The $L_{1}$ rate of contraction of the pseudo posterior distribution with respect to the joint distribution for population generation and the taking of informative samples is derived.  Our approach includes two unique enabling results.
Please see Appendix sections~\ref{AppMain} and \ref{AppEnabling} for details.

\section{Application}\label{application}
We construct a model for count data and perform inference on survey responses collected by the Job Openings and Labor Turnover Survey (JOLTS), introduced in Example~$5$ of Section~\ref{examples}, which is administered by BLS on a monthly basis to a randomly-selected sample from a frame composed of non-agricultural U.S. private (business) and public establishments.  JOLTS focuses on the demand side of U.S. labor force dynamics and measures job hires, separations (e.g. quits, layoffs and discharges) and openings.  The JOLTS sampling design assigns inclusion probabilities (under sampling without replacement) to establishments to be proportional to the number of employees for each establishment (as obtained from the Quarterly Census of Employment and Wages (QCEW)).  This design is informative in that the number of employees for an establishment will generally be correlated with the number of hires, separations and openings.  We perform our modeling analysis on a May, $2012$ data set of $n = 8595$ responding establishments.

We begin by specifying a finite population regression probability model from which we formulate the sampling-weighted pseudo posterior joint distribution that we use to make inference on model parameters from the population generating distribution with only the observed sample of a finite population.  We demonstrate that failing to incorporate sampling weights (e.g. by estimating the posterior distribution defined for the finite population on the observed sample) produces large differences in estimates of parameters.

Our regression model defines a multivariate response as the number of job hires (Hires) for the first response variable and total separations (Seps) as the second response variable.  We construct a single multivariate model (as contrasted with the specification of two univariate models) because these variables of interest tend to be highly correlated such that we expect the regression parameters to express dependence; for example, these two variables are correlated at $60\%$ in our May $2012$ dataset.

We formulate a model for count data that accommodates the high degree of over-dispersion expressed in our establishment-indexed multivariate responses due to the large employment size differences across the establishments.  Were we working with domain-indexed (e.g., by state or county) responses, we may consider to use a Gaussian approximation for the count data likelihood, but such is not appropriate for us due to the presence of many small-sized establishments.  The modeling of count data outcomes is very typical for the analysis of BLS survey data for establishments focused on (un)employment.

We specify the following count data model for the population,
\begin{align}
y_{id} &\ind \mbox{Pois}\left(\exp\left(\psi_{id}\right)\right)\\
\mathop{\Psi}^{N\times D} &\sim \mathop{\mathbf{X}}^{N\times D}\mathop{\mathbf{B}}^{P\times D} + \mathcal{N}_{N\times D}\left(\mathbb{I}_{N},\mathop{\Lambda^{-1}}^{D\times D}\right)\label{likeB}\\
\mathop{\mathbf{B}} &\sim \mathbf{0}
            + \mathcal{N}_{P\times D}\left(\mathop{\mathbf{M}^{-1}}^{P\times P},\left[\tau_{B}\Lambda\right]^{-1}\right)\label{priorB}\\
\Lambda &\sim \mathcal{W}_{D}\left((D+1),\mathbf{I}_{D}\right)\\
\tau_{B} &\sim \mathcal{G}\left(1,1\right)\\
\mathbf{M} &\sim \mathcal{W}_{P}\left((P+1),\mathbf{I}_{P}\right),
\end{align}
where $i = 1,\ldots,N$ indexes the number of establishments and $d = 1,\ldots,D$ indexes the number of dimensions for the multivariate response, $\mathbf{Y}$.  The $N\times D$ log-mean, $\displaystyle\Psi = \left(\mathop{\bm{\psi}_{1}^{'}}^{D\times 1},\ldots,\bm{\psi}_{N}^{'}\right)$, may be viewed as a latent response whose columns index the number of job hires (Hires) and total separations (Seps) under our JOLTS application, so that $D=2$.  The number of predictors in the design matrix, $\mathbf{X}$, is denoted by $P$ and $\mathbf{B}$ are the unknown matrix of population coefficients that serve as the focus for our inference.  Our model is formulated as a multivariate Poisson-lognormal model, under which the Gaussian prior of Equation~\ref{likeB} for the logarithm of the Poisson mean allows for over-dispersion (of different degrees) in each of the $D$ dimensions.  The priors in Equation~\ref{likeB} and Equation~\ref{priorB} are formulated in matrix variate (or, more generally, tensor product) Gaussian distributions using the notation of \citet{dawid:1981}; for example, the prior for the $P \times D$ matrix of coefficients, $\mathbf{B}$, assigns the $P\times D$ mean $\mathbf{0}$ for a Gaussian distribution that employs a separable covariance structure where the $P \times P,~\mathbf{M}$, denotes the precision matrix for the columns of $\mathbf{B}$, and the $D\times D,~\tau_{B}\Lambda$, denotes the precision matrix for the rows.  This prior formulation is the equivalent of assigning a $PD$ dimensional Gaussian prior to a vectorization of $\mathbf{B}$ accomplished by stacking its columns with $PD\times PD$ precision matrix, $\mathbf{M}\otimes\left(\tau_{B}\Lambda\right)$. (See \citet{hoff2011} for more background).  Precision matrices, $\left(\mathbf{M},\Lambda\right)$, each receive Wishart priors with hyperparameter values that impose uniform marginal prior distributions on the correlations \citep{CIS-168884}.

We regress the multivariate latent response, $\Psi$, on predictors representing the logarithm of the overall establishment-indexed number of employees (Emp), obtained from the QCEW, the logarithm of the number of job openings (Open), region (Northeast, South, West, Midwest (Midw)) and ownership type (Private, Federal Government, State Government (State), Local Government (Local)).  We convert region and ownership type to binary indicators and leave out the Northeast region and Federal Government ownership to provide the baseline of a full-column rank predictor matrix.  We summarize our regression model on the logarithm scale by: $(\psi_{\mbox{\tiny{Hires}}}, \psi_{\mbox{\tiny{Seps}}})$ $\sim$  $1$ + West $+$ Midw $+$ South $+$ State $+$ Local $+$ Private $+$ $\log(\mbox{Emp})$ $+$ $\log(\mbox{Opens})$ $+$ error, where $1$ denotes an intercept (Int).

Our population model is hypothesized to generate the finite population of the U.S. non-agricultural establishments, from which we have taken a sample of size $n = 8595$ for May, $2012$ as our observations.  For ease of reading, we will continue to use $\mathbf{Y}$ and $\mathbf{X}$, to next define the associated pseudo posterior, though each possesses $n < N$ rows representing the sampled observations, in this context.

The population model likelihood contribution for establishment, $i$, on dimension, $d$, is formed with the integration,
\begin{equation}
p\left(y_{id}\middle\vert\mathbf{x}_{i},\mathbf{B},\Lambda\right) = \mathop{\int}_{\mathbb{R}}p\left(y_{id}\middle\vert\psi_{id}\right)\times p\left(\psi_{id}\middle\vert\mathbf{x}_{i},\mathbf{B},\Lambda\right)d\psi_{id},
\end{equation}
where sampling weight, $w_{i} = 1/\pi_{i}$ and $\tilde{w}_{i} = n \times w_{i}/\sum_{i=1}^{n}w_{i}$, such that the adjusted weights sum to $n$, the asymptotic amount of information contained in the sample (under a sampling design that obeys condition~\nameref{independence}).  This integrated likelihood induces the following pseudo likelihood,
\begin{equation}
p^{\pi}\left(y_{id}\middle\vert\mathbf{x}_{i},\mathbf{B},\Lambda\right) = \left[\mathop{\int}_{\mathbb{R}}p\left(y_{id}\middle\vert\psi_{id}\right)\times p\left(\psi_{id}\middle\vert\mathbf{x}_{i},\mathbf{B},\Lambda\right)d\psi_{id}\right]^{\tilde{w}_{i}},
\end{equation}
which is analytically intractable, so we perform the integration, numerically, in our MCMC using the prior
for each $\psi_{id}$ exponentiated by the normalized sampling weight, $\tilde{w}_{i}$, which we use to construct its pseudo posterior distribution.  Using Bayes rule we present the logarithm of the pseudo posteriors for the latent set of $D\times 1$ log-mean parameters, $\{\bm{\psi}_{i}\}$, (which are \emph{a posteriori} independent over $i = 1,\ldots,n$), with,
\begin{subequations}
\begin{align}
&\log~p^{\pi}\left(\bm{\psi}_{i}\middle\vert\mathbf{y}_{i},\mathbf{x}_{i},\mathbf{B},\Lambda\right) \propto\\
&\log\left\{\left[\mathop{\prod}_{d=1}^{D}\exp\left(\psi_{id}\right)^{y_{id}}\exp\left(-\exp\left(\psi_{id}\right)\right)\right]^{\tilde{w}_{i}}
\times \left[\mathcal{N}_{D}\left(\bm{\psi}_{i}\middle\vert\mathbf{x}_{i}^{'}\mathbf{B},\Lambda^{-1}\right)\right]^{\tilde{w}_{i}}\right\}\\
&\propto\tilde{w}_{i}\sum_{d=1}^{D}\left[y_{id}\psi_{id} - \exp\left(\psi_{id}\right)\right] - \frac{1}{2}
\left(\bm{\psi}_{i}-\mathbf{x}_{i}^{'}\mathbf{B}\right)^{'}\tilde{w}_{i}\Lambda\left(\bm{\psi}_{i}-\mathbf{x}_{i}^{'}\mathbf{B}\right)
\label{loglikepsi},
\end{align}
\end{subequations}
where we note in the second expression in Equation~\ref{loglikepsi} that the sampling weights influence the \emph{prior} precision for each $\bm{\psi}_{i}$, such that a higher-weighted observation will exert relatively more influence on posterior inference because this observation is relatively more representative of the population.  We take samples from the pseudo posterior distribution specified Equation~\ref{loglikepsi} in our MCMC using the elliptical slice sampler of \citet{murray2010}, where we draw $\bm{\psi}_{i} \ind \mathcal{N}_{D}\left(\mathbf{x}_{i}^{'}\mathbf{B},\left(\tilde{w}_{i}\Lambda\right)^{-1}\right)$ and formulate a proposal as a convex combination (parameterized on an ellipse) of this draw from the prior and the value selected on the previous iteration of the MCMC.   We evaluate each proposal using the weighted likelihood in the first expression of Equation~\ref{loglikepsi}.

We next illustrate the construction of the pseudo posterior distribution for the $P\times D$ matrix of regression coefficients, $\mathbf{B}$, (which by D-separation is independent of the observations, ($(y_{id})$, given $(\psi_{id})$),
\begin{subequations}
\begin{align}
p^{\pi}\left(\mathbf{B}\vert \mathbf{Y},\mathbf{X},\Psi,\Lambda,\mathbf{M},\tau_{B}\right) &\propto
    \left[\mathop{\prod}_{i=1}^{n}\mathcal{N}_{n\times D}\left(\bm{\psi}_{i}\vert\mathbf{B}^{'}\mathbf{x}_{i},\mathbb{I}_{n},\Lambda^{-1}\right)^{\tilde{w}_{i}}\right]
    \mathcal{N}_{P\times D}\left(\mathbf{B}\vert \mathbf{M}^{-1},\left(\tau_{B}\Lambda\right)^{-1}\right)\\
\log~p^{\pi}\left(\mathbf{B}\vert \mathbf{Y},\mathbf{X},\Psi,\Lambda,\mathbf{M},\tau_{B}\right) &\propto \mathop{\sum}_{i=1}^{n}\left[\frac{\tilde{w}_{i}}{2}\log\vert\Lambda\vert - \frac{\tilde{w}_{i}}{2}\left(\bm{\psi}_{i}-\mathbf{B}^{'}\mathbf{x_i}\right)^{'}\Lambda
    \left(\bm{\psi}_{i}-\mathbf{B}^{'}\mathbf{x_i}\right)\right] \nonumber\\
    &+ \log~\mathcal{N}_{P\times D}\left(\mathbf{B}\vert \mathbf{M}^{-1},\left(\tau_{B}\Lambda\right)^{-1}\right)\label{loglikeB}.
\end{align}
\end{subequations}
In a Bayesian setting, the sum of the weights ($n = \sum_{i=1}^{n}\tilde{w}_{i}$) impacts the estimated posterior variance as we observe in Equation~\ref{loglikeB}.  We see that weights scale the quadratic product of the Gaussian kernel in Equation~\ref{loglikeB} such that we may accomplish the same result using the matrix variate formation to define the pseudo likelihood, $\mathcal{N}_{n\times D}\left(\Psi-\mathbf{X}\mathbf{B}\vert\mathbf{\tilde{W}},\Lambda^{-1}\right)$, where $\mathbf{\tilde{W}}= \mbox{diag}\left(\tilde{w}_{1},\ldots,\tilde{w}_{n}\right)$, the weights for the sampled observations, from which we compute the following conjugate conditional pseudo posterior distribution defined on the $n$ observations,
\begin{equation}\label{pseudoB}
p^{\pi}\left(\mathbf{B}\vert\mathbf{Y},\mathbf{X},\Psi,\Lambda,\mathbf{M},\tau_{B}\right) = \mathbf{h}^{\pi}_{B} +
    \mathcal{N}_{P\times D}\left(\mathbf{B}\vert (\bm{\phi}^{\pi}_{B})^{-1},\Lambda^{-1}\right),
\end{equation}
where $\bm{\phi}^{\pi}_{B} = \mathbf{X}^{'}\mathbf{\tilde{W}}\mathbf{X} + \tau_{B}\mathbf{M}$ and $\mathbf{h}^{\pi}_{B} =
(\bm{\phi}^{\pi}_{B})^{-1}\mathbf{X}^{'}\mathbf{\tilde{W}}\Psi$.

Under employment of a simpler continuous response framework, the conditional posterior for $\mathbf{B}$ retains the same form as Equation~\ref{pseudoB}, except the latent response on the logarithm scale, $\Psi$, would be replaced by the observed data, $\mathbf{Y}$.  Intuitively, we note using a sampling-weighted pseudo prior for the latent response, $\Psi$, for sampling coefficients, $\mathbf{B}$, is analogous to using the sampling-weighted likelihood in the case of an observed, continuous response, $\mathbf{Y}$.

\begin{figure}[!Ht]
\begin{center}
\includegraphics[width=5.4in,height=3.8in]{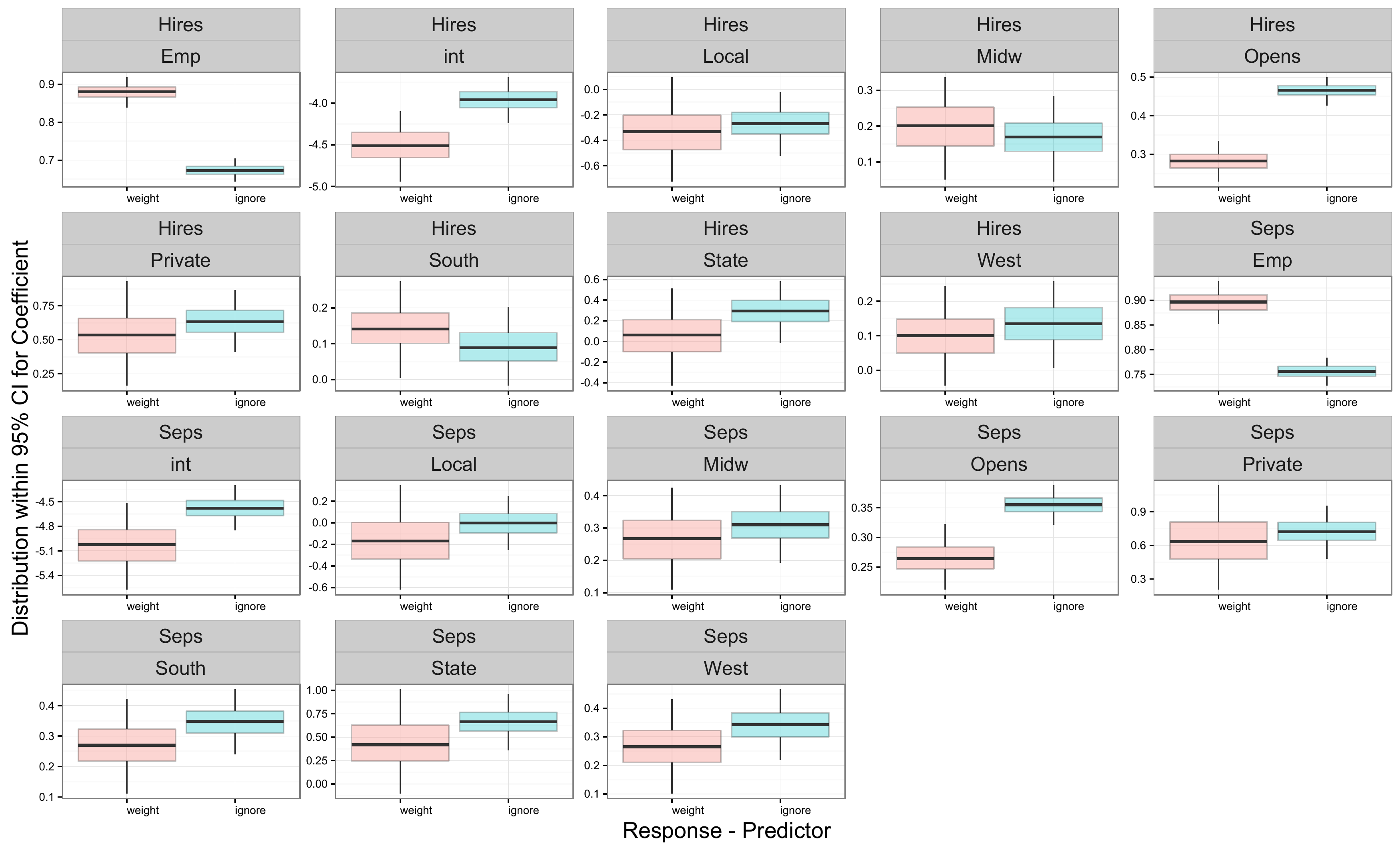}
\caption{Comparison of posterior densities for the each coefficient in the $(P=9)\times (D=2)$ coefficient matrix,
$\mathbf{B}$, within $95\%$ credible intervals, based on inclusion sampling weights in a pseudo posterior (the left-hand plot in each panel) and exclusion
of the sampling weights using the posterior distribution defined for the population (in the right-hand plot).
Each plot panel is labeled by ``predictor,response" for the two included response variables, ``Hires", and ``Seps" (total separations).}
\label{joltsresults}
\end{center}
\end{figure}

Each plot panel in Figure~\ref{joltsresults} compares estimated posterior distributions for a coefficient in $\mathbf{B}$ (within $95\%$ credible intervals), labeled by ``predictor, dimension (of the multivariate response)", when applied to the May, $2012$ JOLTS dataset between two estimation models: 1. The left-hand plot in each panel employs the sampling weights to estimate the pseudo posterior for $\mathbf{B}$, induced by the pseudo posterior for the latent response in Equation~\ref{loglikepsi}; 2. The right-hand plot estimates the coefficients using the posterior distribution defined on the finite population, which may be achieved by replacing $\mathbf{\tilde{W}}$ by the identity matrix to equally weight establishments.  Equal weighting of establishments assumes that the sample represents the same balance of information as the population from which it was drawn, which is not the case under an informative sampling design.  Comparing estimation results from the pseudo posterior and population posterior distributions provides one method to assess the sensitivity of estimated parameter distributions to the sampling design.

We observe that the estimated results are quite different in both location and variation between estimations performed under the pseudo posterior and population posterior distributions, indicating a high degree of informativeness in the sampling design.  The $95\%$ credible intervals for the coefficients of the continuous predictors - (the log of) job openings (Opens) and employment (Emp) - don't even overlap on both the number of hires (Hires) and separations (Seps) responses.  The coefficient for the State ownership predictor and the number of hires response is bounded away from $0$ when estimated under the (unweighted) population posterior, but is centered on $0$ under the sampling-weighted, pseudo posterior.  The coefficient posterior variances estimated on the observed sample under the population posterior are understated because they don't reflect the uncertainty with which the information in the sample expresses that in the population (which is captured through the sampling weights).
\FloatBarrier

\subsection{Simulation Study}\label{simulation}

We implement a simulation study to compare the marginal \emph{pseudo} posterior distributions to the (unweighted) population posterior distributions for the regression coefficients, where both are estimated on the observed sample drawn under an informative sampling design.
For this study we use the $N = 8595$ observations from the JOLTS May, $2012$ data as our \emph{population}.
We take $100$ Monte Carlo samples of size $\mathbf{n}_{\nu} = \left(500,1000,1500,2500\right)$ establishments
using an informative single-stage sample design with unequal inclusion probabilities based on the proportional to size sample used for the real JOLTS survey.
Characteristics of the the sampling design, used for this study, at each sample size are presented in Table~\ref{tab:design}.

\begin{figure}[!Ht]
\begin{center}
\includegraphics[width=5.0in,height=3.4in]{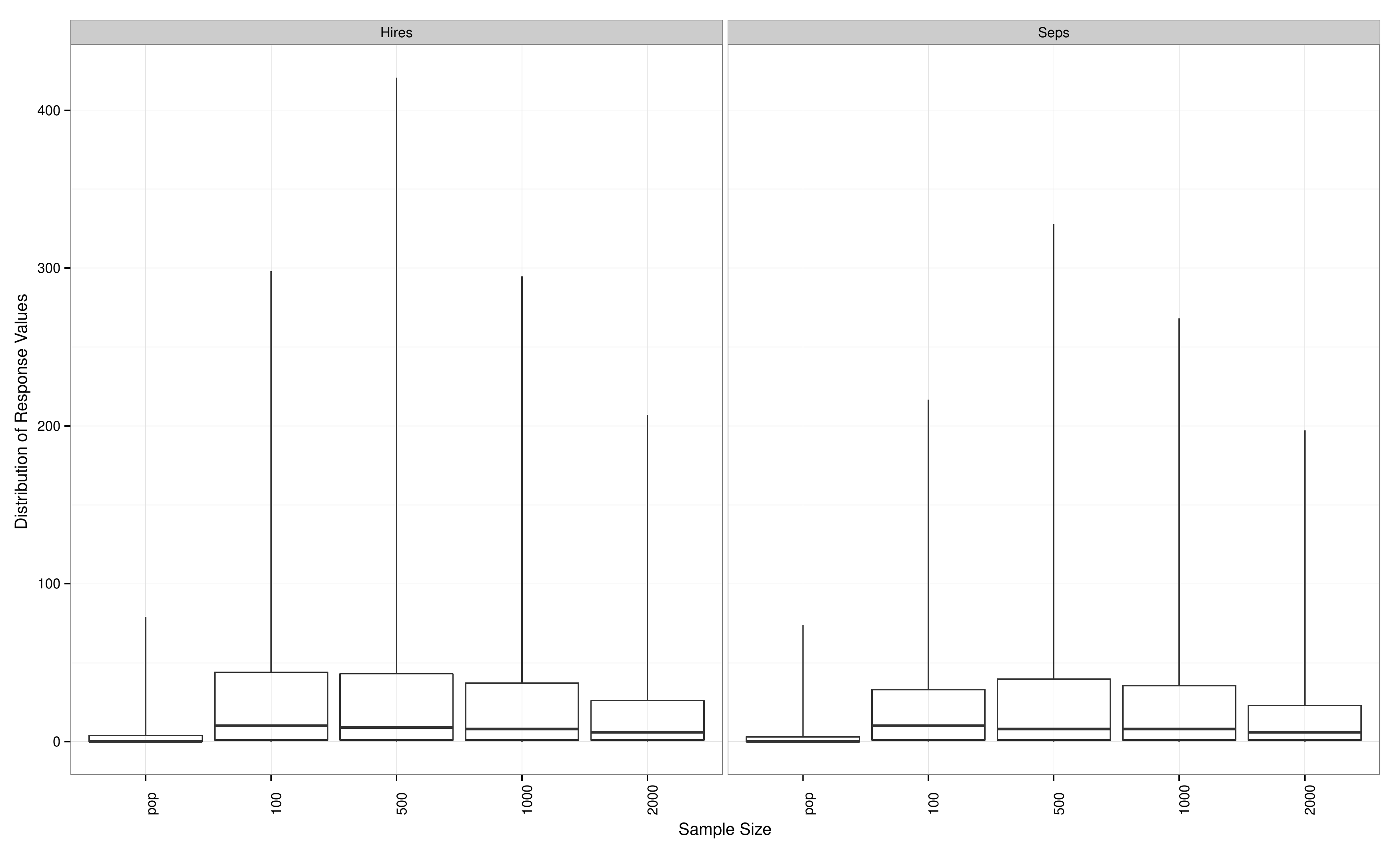}
\caption{Distributions of response values for population compared to informative samples.  The left-most box plot in each of the two plot panels contains the distribution for the JOLTS sample that we use as our ``population" in the simulation study.   The next set of box plots show the distribution for the response values for increasing sample sizes (from left-to-right) for each sample drawn under our single stage proportion-to-size design.  The left-hand plot panel displays the Hires response variable and the right-hand panel displays the Seps (separations) response variable.}
\label{infbal}
\end{center}
\end{figure}

This sampling design will induce distributions of the observed samples that will be different from those for the population.  The designed correlation between the response and inclusion probabilities will produce observed samples with values skewed towards higher numbers of hires and separations than in the population.  Figure~\ref{infbal} demonstrates this difference between the distributions for realized samples under the informative sampling design compared to those for the finite population.  The left-most box plot in each of the two panels displays the population distribution for a response value.  A single sample is drawn under a sequence of increasing sample sizes for illustration.  The next set of box plots displays the resulting distributions for the response values in each sample with size increasing from left-to-right.  The left-hand plot panel displays the distributions for the Hires response, while the right-hand panel displays those for the Seps (separations) response variable.

\begin{figure}[!Ht]
\begin{center}
\includegraphics[width=5.4in,height=3.8in]{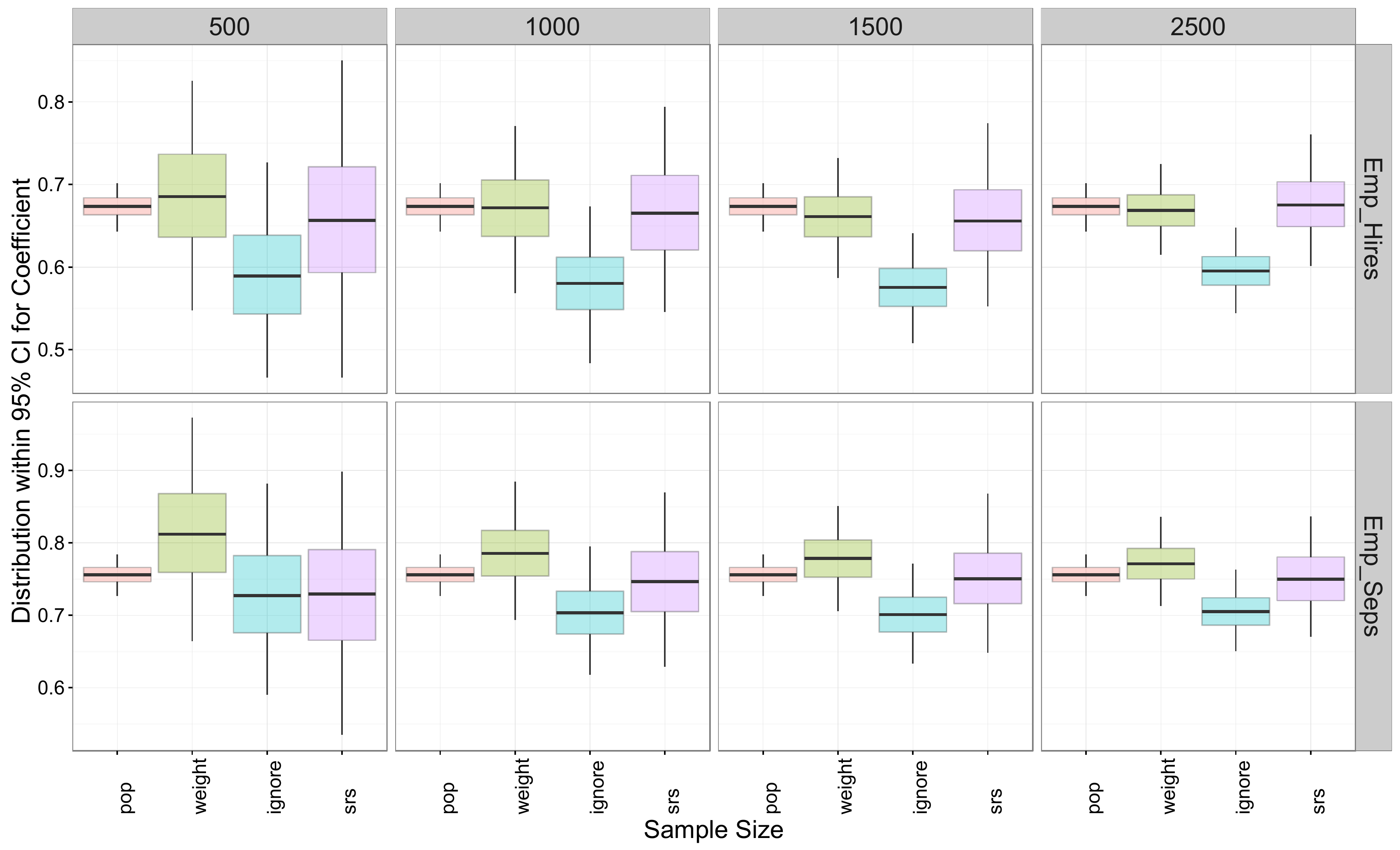}
\caption{Comparison of posterior densities for $2$ coefficients, Employment-Hires (top row of plot panels) and Employment-Separation (bottom row of plot panels) in $\mathbf{B}$, within $95\%$ credible intervals, between estimation on the population (left-hand plot in each panel), estimations from informative samples data taken from that population, which include sampling weights in a pseudo posterior (the second plot from the left in each panel) and exclusion
of the sampling weights using the population posterior distribution (the third plot from the left) under a simulation study. The right-most plot presents the posterior density estimated from a simple random sample of the same size for comparison.  The simulation study uses the May, $2012$ JOLTS sample as the ``population" and generates $500$ informative samples for a range of sample sizes (of $500, 1000, 1500, 2500$, from left-to-right) under a sampling without replacement design with inclusion probabilities set proportionally to the square root of employment levels.  A separate estimation is performed on each Monte Carlo sample and the draws from estimated distributions are concatenated over the samples.}
\label{mcsim}
\end{center}
\end{figure}

Pseudo posterior and population posterior distributions are estimated on each Monte Carlo sample at each sample size in $\mathbf{n}_{\nu}$. Figure~\ref{mcsim} compares estimation of the posterior distribution from the fully-observed population (left-hand box plot) to estimation using the \emph{pseudo} posterior from sample observations taken under the proportional-to-size sampling design. The third box plot in each panel shows the estimation of the posterior distribution estimated on the same sample ignoring the informative sampling design.  The last box plot in each panel displays the estimates of the posterior distribution from a simple random sample of the same size, where no correction for the sampling design is required,  as a gold standard against which to measure the performance of the pseudo posterior distribution.  We estimate the distributions on each of the $100$ Monte Carlo draws for each sample size and concatenate the results such that they incorporate both the variation of population generation and repeated sampling from that population.  The sample sizes, $n_{\nu}$, increase from left-to-right across the plot panels.  The top set of plot panels display the posterior distributions of the regression coefficient for the employment predictor (Emp) and the hires response (Hires), while the bottom set of panels display the coefficient distributions for the employment predictor (Emp) and the total separations response (Seps).

\begin{table}[!Ht]
\centering
\begin{tabular}{rrrrrrrr}
  \hline
  \hline
 & $n_{\nu}$ & CUs & min($\pi_{\nu}$) & max($\pi_{\nu}$) & $\mbox{CV}(\pi_{\nu})$ & Cor($y_{\mbox{\tiny{hires}}},\pi_{\nu}$) & Cor($y_{\mbox{\tiny{Seps}}},\pi_{\nu}$) \\
  \hline
1 & 500 &  56 & 0.02 & 1.00 & 2.11 & 0.80 & 0.62 \\
  2 & 1000 & 196 & 0.04 & 1.00 & 1.60 & 0.69 & 0.50 \\
  3 & 1500 & 357 & 0.07 & 1.00 & 1.29 & 0.61 & 0.44 \\
  4 & 2500 & 722 & 0.14 & 1.00 & 0.91 & 0.51 & 0.36 \\
   \hline
\end{tabular}
\caption{Characteristics of single stage, fixed size pps sampling design used in simulation study.  $n_{\nu}$ denotes
the sample size. CUs denotes the number of certainty units (with inclusion probabilities equal to 1). $\pi_{\nu}$ denotes the inclusion probabilities (proportional to square root of JOLTS employment), $\mbox{CV}(\pi_{\nu})$ denotes the coefficient of variation of $\pi_{\nu}$, Cor($y_{\mbox{\tiny{hires}}}$,$\pi_{\nu}$) denotes correlation of the number of hires and $\pi_{\nu}$ and Cor($y_{\mbox{\tiny{Seps}}}$,$\pi_{\nu}$) denotes the correlation of the number of separations
and $\pi_{\nu}$.}
\end{table}\label{tab:design}

\FloatBarrier

\FloatBarrier

Scanning from left-to-right in each row across the increasing sample sizes, we readily note a consistent difference in the estimated posterior mean, as expected, between the population model estimated on the samples without adjustment for the informative sampling design as compared to the mean of the posterior distribution estimated on the entire finite population.  The application of the pseudo posterior model, however, produces much less difference (relative to estimation on the fully observed population), though the difference between the estimated pseudo-posterior and the population posterior is yet notably more than that for the simple random sampling result (estimated on samples of the same size as the informative sample).  The estimated difference for the pseudo-posterior converges to $0$, however, as the sample size increases.  The posterior variance for the estimated posterior under simple random sampling remains larger than that for the pseudo posterior estimated on the informative sample because our proportion-to-size sampling design over-samples the highest variance units, which provides better capture of information in the population (which is why this design is used).  So, in this case, the improved capture of information in the finite population provided by our sampling design more than overcomes the added variation induced by estimation with the sampling weights.  In summary, this simulation study demonstrates the contraction of the pseudo posterior distribution estimated on the sample onto the posterior distribution estimated on a fully-observed finite population.

We were able to directly perform posterior inference about the population using only quantities available for the observed sample under the pseudo posterior full conditional distributions outlined in Section~\ref{application}.  By contrast, \citet{Litt:mode:2004} offer no modeling approach that parameterizes a proportion-to-size sampling design because they note that each unit is in its own group under the stratum-indexed construction they generally suggest.  A typical naive approach, however, is to simply include the sampling weights or a variable highly corrected with them as a predictor \emph{only} for observed units with no imputation of the non-sampled units.  This is precisely the construction of the alternative that estimates the population posterior distribution on the informative sample, which is shown as the third box plot in each plot panel of Figure~\ref{mcsim}, because the employment variable, Emp, is included as a predictor and is highly correlated with the sampling weights.  This option includes Emp only for the sampled units as does the model for the pseudo posterior.   The reason for biased inference, even when including a predictor that is highly correlated with sampling weights, is because the distribution for the sampled data conditioned on the sampling weights is not generally equal to the distribution for the population conditioned on the sampling weights by Bayes rule \citep{Pfef:Sver:infe:2009}.

The JOLTS respondent-level data from which samples were drawn for our Monte Carlo simulation study may not be publicly released due to restrictions that protect confidentiality of survey participants.  A Monte Carlo simulation study using our pseudo posterior plug-in method may, however, be generated under a Bayesian nonparametric model for functional data that is available from the \texttt{growfunctions} package for R \citep{savitskyr:2016}. The package includes a synthetic data generator whose use is illustrated in \citet{savitskyjss:2014}, along with a Monte Carlo simulator that compares parameter estimates when correcting versus ignoring an informative sampling design.  Although the functional data model is more complicated than our count data application, the Monte Carlo simulation function available in \texttt{growfunctions} produces a figure that demonstrates results very similar to those displayed in Figure~\ref{mcsim}.

\section{Discussion}\label{discussion}
A variety of broadly applicable approaches are available for incorporating sampling weights into weighted maximum likelihood estimation procedures \citep{Pfef:Sver:infe:2009} to account for an informative sampling design.  Defining easily adaptable algorithms that account for sampling design informativeness under any model for the population specified by the data analyst has proved more challenging for estimating Bayesian probability models.  Solutions have focused on parameterizing the sampling design or co-estimating the conditional expectation of inclusion, along with the population-generating model.  While these approaches allow estimation using the sampled observations, the implementations typically require a high degree of customization to each population model. We take a different approach that constructs a  sample-weighted pseudo-posterior to account for an informative sampling design in our plug-in method that is readily accommodated to any Bayesian population probability model.

We demonstrated the applicability of the plug-in method to a poisson - lognormal model for count data.  We showed that the plug-in method reduces estimation bias and posterior estimation includes the uncertainty with which the sample reflects the population on these covariance parameters.

The plug-in method is as easily-implemented and broadly applicable as those methods used for likelihood based optimization.  We illustrated in Section~\ref{application} that the full conditional posterior distributions defined for the population generating model are easily updated by multiplying the log-likelihoods for $\left(\{y_{id}\},\psi_{id}\right)$, by the sampling weights, $\{\tilde{w}_{i}\}$, without changing the constructions for full conditional posterior distributions.  The same concerns that apply in the use of sampling weights under likelihood optimization also apply for Bayesian estimation. The quality of posterior estimation is highly dependent on the population representativeness of the realized sample.  Sampling weights may be adjusted based on the composition of the realized sample through estimating the conditional expectation of the weights, given the response values for the observed units.  Regressing the weights on the response variables using the observed data and replacing the raw weights with their conditional expectation, known as ``weight smoothing", would be expected to reduce the posterior variances for estimated parameters to the extent that the weights express variance unrelated to the response.  While the conditional distribution for the sampling weights given the response under the realized sample is not generally expected to be the same as that for the finite population, their expectations are equal \citep{Pfef:Sver:infe:2009}. We explored such smoothing for our sampling weights for the JOLTS application, but there was little reduction in variance, so we employed the published weights for simplicity.

Even after adjustment, if the composition of the realized sample unevenly reflects information in the population, the weights would express a high variation. Approaches that calibrate the sampling weights to actual population totals, where known, may improve the quality of estimation produced from the plug-in method.  BLS performs a calibration adjustment of the JOLTS sampling weights such that the weighted difference of hires and separations reported in the sample ties to the monthly total employment change from the CES survey. (The CES survey is introduced and discussed in Section~\ref{examples}).  This step adjusts the sampling weights computed from inclusion probabilities under the JOLTS sampling design based on the \emph{actual} sample achieved in each month.

One may, alternatively, implement a more fully Bayesian approach that parameterizes a joint model for the response and sampling weights, specific to a given population generating model, as a method that smoothes the weights in the presence of the response values.  Doing so, however, requires imputation of weights and response values for non-sampled units, which can be computationally expensive for a survey that samples from the entire U.S. population of business establishments, as does JOLTS.

Lastly, we construct conditions which, together, define a class of sampling designs under which $L_{1}$ consistency of the pseudo posterior is guaranteed.  One of these conditions requires that the pairwise sample inclusion dependencies asymptotically decrease to $0$.  While there are many sampling designs, in practice, which are members of this class, including the proportion-to-size sampling design used for our JOLTS application, there are some designs which are not - such as a cluster sampling design where the number of clusters grows, but the number of units in each cluster remains relatively fixed.   A direction for future research will be to widen the class of allowable designs by incorporating second order (or pairwise) inclusion probabilities for inference, though doing so will introduce some practical restrictions on the specifications for the population generating model.

\bibliography{mv_refs_sep2015}
\bibliographystyle{agsm}

\newpage
\appendix
\section{Proof of Theorem~\ref{results}}\label{AppMain}
\begin{proof}
Condition~\nameref{existtests} establishes the existence of test statistics, $\phi_{n_{\nu}}\left(X_{1}\delta_{\nu 1},\ldots,X_{N_{\nu}}\delta_{\nu N_{\nu}}\right) \in \left(0,1\right)$ used to achieve the
following result,
\begin{align}
&\mathbb{E}_{P_{0},P_{\nu}}\phi_{n_{\nu}} \nonumber\\
&\leq\exp\left(n_{\nu}\xi_{N_{\nu}}^{2}\right)\cdot\frac{\exp\left(-K n_{\nu}M^{2}\xi_{n_{\nu}}^2\right)}{1-\exp\left(-K n_{\nu}M^{2}\xi_{N_{\nu}}^2\right)}\nonumber\\
&\leq 2\exp\left(-K n_{\nu}\xi_{N_{\nu}}^{2}\right)\label{exptest},
\end{align}
in Lemmas $2$ and $9$ of \citet{ghosal2007} by setting $\xi = M\xi_{N_{\nu}}$, and by choosing constant $M >0$ sufficiently large, such that $KM^{2}-1> K$.

We will bound the expectation (under $\left(P_{0},P_{\nu}\right)$, jointly) of the mass assigned by pseudo posterior distribution for those $P$ at some minimum distance from $P_{0}$,
\begin{align}
&\displaystyle\Pi^{\pi}\left(P\in\mathcal{P}:d^{\pi}_{N_{\nu}}\left(P,P_{0}\right) \geq M\xi_{N_{\nu}} \vert X_{1}\delta_{\nu 1},\ldots,X_{N_{\nu}}\delta_{\nu N_{\nu}}\right)\nonumber \\
&= \Pi^{\pi}\left(P\in\mathcal{P}:d^{\pi}_{N_{\nu}}\left(P,P_{0}\right) \geq M\xi_{N_{\nu}} \vert X_{1}\delta_{\nu 1},\ldots,X_{N_{\nu}}\delta_{\nu N_{\nu}}\right)\left(\phi_{n_{\nu}} + 1-\phi_{n_{\nu}}\right)\label{totalpost}.
\end{align}

Equation~\ref{exptest} establishes the bound,
\begin{multline}
\mathbb{E}_{P_{0},P_{\nu}}\Pi^{\pi}\left(P\in\mathcal{P}:d^{\pi}_{N_{\nu}}\left(P,P_{0}\right)\geq M\xi_{N_{\nu}}\middle\vert X_{1}\delta_{\nu 1},\ldots,X_{N_{\nu}}\delta_{N_{\nu}}\right)\phi_{n_{\nu}} \leq \mathbb{E}_{P_{0}}\phi_{n_{\nu}}\leq\\
2\exp\left(-K n_{\nu}\xi_{N_{\nu}}^{2}\right),
\end{multline}
since the pseudo posterior mass is bounded from above by $1$.   We next enumerate the pseudo posterior distribution for the second term of Equation~\ref{totalpost},
\begin{multline}\label{post}
\Pi^{\pi}\left(P\in\mathcal{P}:d^{\pi}_{N_{\nu}}\left(P,P_{0}\right)\geq M\xi_{N_{\nu}}\middle\vert\mathbf{X}_{N_{\nu}}\bm{\delta}_{N_{\nu}}\right)\left(1-\phi_{n_{\nu}}\right)=\\
\frac{\displaystyle\mathop{\int}_{P\in\mathcal{P}:d^{\pi}_{N_{\nu}}\left(P,P_{0}\right)\geq M\xi_{N_{\nu}}}\mathop{\prod}_{i=1}^{N_{\nu}}\frac{p^{\pi}}{p_{0}^{\pi}}\left(X_{i}\delta_{\nu i}\right)d\Pi\left(P\right)\left(1-\phi_{n_{\nu}}\right)}
{\displaystyle\mathop{\int}_{P\in\mathcal{P}}\mathop{\prod}_{i=1}^{N_{\nu}}\frac{p^{\pi}}{p_{0}^{\pi}}
\left(X_{i}\delta_{\nu i}\right)d\Pi\left(P\right)}.
\end{multline}

We may bound the denominator of Equation~\ref{post} from below, in probability.  Define the event,
\begin{equation*}
B_{N_{\nu}} = \left\{P:-P_{0}\log\left(\frac{p}{p_{0}}\right) \leq \xi_{N_{\nu}}^2, P_{0}\left(\log\frac{p}{p_{0}}\right)^{2} \leq \xi_{N_{\nu}}^{2}\right\}
\end{equation*}

We have from Lemma~\ref{denominator},
\begin{equation*}
\mbox{Pr}\left\{\mathop{\int}_{P\in\mathcal{P}}\displaystyle\mathop{\prod}_{i=1}^{N_{\nu}}\frac{p^{\pi}}{p_{0}^{\pi}}
\left(X_{i}\delta_{\nu i}\right)d\Pi\left(P\right)\geq \exp\left[-(1+C)N_{\nu}\xi^{2}\right]\right\}
\geq 1 - \frac{\gamma+C_{3}}{C ^{2}N_{\nu}\xi^{2}},
\end{equation*}
for every $P \in B_{N_{\nu}}$ and any $C > 0$, $\gamma > 1$, where $\gamma$ may be set closer to $1$ for sampling designs that define a low gradient for inclusion probabilities, $\{\pi_{\nu i}\}$.  The constant, $C_{3} > 0$, and will be close to $1$ for sufficiently large $\nu$.  Condition~\nameref{priortruth} restricts the prior on $B_{N_{\nu}}$,
\begin{equation*}
\Pi\left(B_{N_{\nu}}\right) \geq \exp\left(-N_{\nu}\xi_{N_{\nu}}^{2}C\right).
\end{equation*}
Then with probability at least $1 - \frac{16\gamma^{2}\left[\gamma+C_{3}\right]}{\left(K M^{2}f - 2\gamma\right)^{2} N_{\nu}\xi^{2}}$,
\begin{align*}
\displaystyle\mathop{\int}_{P\in\mathcal{P}}\mathop{\prod}_{i=1}^{N_{\nu}}\frac{p^{\pi}}{p_{0}^{\pi}}
\left(X_{i}\delta_{\nu i}\right)d\Pi\left(P\right)
&\geq\exp\left[-(1+C)N_{\nu}\xi^{2}\right]\Pi\left(B_{N_{\nu}}\right)\\
& \geq \exp\left(-(1+2C)N_{\nu}\xi^{2}\right)\\
& \geq \exp\left(-\frac{K M^{2} n_{\nu}\xi_{N_{\nu}}^{2}}{2\gamma}\right),
\end{align*}
where we set $1+2C = \frac{K M^{2}f}{2\gamma}$, where we use condition~\nameref{fraction} to replace $f\times N_{\nu}$ with ${n_{\nu}}$ for $\nu$ sufficiently large.

Denote this event by,
\begin{equation}
A^{\pi}_{N_{\nu}} = \left\{\displaystyle\mathop{\int}_{P\in\mathcal{P}}\mathop{\prod}_{i=1}^{N_{\nu}}\frac{p^{\pi}}{p_{0}^{\pi}}
\left(X_{i}\delta_{\nu i}\right)d\Pi\left(P\right)\geq \exp\left(-\frac{K M^{2} n_{\nu}\xi_{N_{\nu}}^{2}}{2\gamma}\right)\right\},
\end{equation}
such that,
\begin{align*}
&\Pi^{\pi}\left(P\in\mathcal{P}:d^{\pi}_{N_{\nu}}\left(P,P_{0}\right)\geq M\xi_{N_{\nu}}\middle\vert\mathbf{X}_{N_{\nu}}\bm{\delta}_{N_{\nu}}\right)\left(1-\phi_{n_{\nu}}\right)\nonumber\\
&=\displaystyle\left[\frac{\displaystyle\mathop{\int}_{\left\{P\in\mathcal{P}:d^{\pi}_{n_{\nu}}\left(P,P_{0}\right)\geq M\xi_{N_{\nu}}\right\}}\mathop{\prod}_{i=1}^{N_{\nu}}\frac{p^{\pi}}{p_{0}^{\pi}}\left(X_{i}\delta_{\nu i}\right)d\Pi\left(P\right)\left(1-\phi_{N_{\nu}}\right)}
{\displaystyle\mathop{\int}_{P\in\mathcal{P}}\mathop{\prod}_{i=1}^{N_{\nu}}\frac{p^{\pi}}{p_{0}^{\pi}}
\left(X_{i}\delta_{\nu i}\right)d\Pi\left(P\right)}\left(\mathbb{I}\left(\left[A^{\pi}_{N_{\nu}}\right]^{c}\right)
+\mathbb{I}\left(A^{\pi}_{N_{\nu}}\right)\right)\right]\nonumber\\
&\leq \mathbb{I}\left(\left[A^{\pi}_{N_{\nu}}\right]^{c}\right) + \mathbb{I}\left(A^{\pi}_{N_{\nu}}\right)\times\biggl[
\exp\left(\frac{K M^{2} n_{\nu}\xi_{N_{\nu}}^{2}}{2\gamma}\right)\Pi\left(\mathcal{P}\backslash\mathcal{P}_{N_{\nu}}\right)\\
&+ \exp\left(\frac{K M^{2} n_{\nu}\xi_{N_{\nu}}^{2}}{2\gamma}\right)\displaystyle\mathop{\int}_{\left\{P\in\mathcal{P}_{N_{\nu}}:d^{\pi}_{N_{\nu}}\left(P,P_{0}\right)\geq M\xi_{N_{\nu}}\right\}}\mathop{\prod}_{i=1}^{N_{\nu}}\frac{p^{\pi}}{p_{0}^{\pi}}\left(X_{i}\delta_{\nu i}\right)d\Pi\left(P\right)\left(1-\phi_{n_{\nu}}\right)\biggr]
\end{align*}

Taking the expectation of both sides with respect to the joint distribution, $\left(P_{0},P_{\nu}\right)$,
\begin{align}
&\mathbb{E}_{P_{0},P_{\nu}}\Pi^{\pi}\left(P\in\mathcal{P}:d^{\pi}_{N_{\nu}}\left(P,P_{0}\right)\geq M\xi_{N_{\nu}}\middle\vert\mathbf{X}_{N_{\nu}}\bm{\delta}_{N_{\nu}}\right)\left(1-\phi_{n_{\nu}}\right)\nonumber\\
&\leq P\left(\left[A^{\pi}_{N_{\nu}}\right]^{c}\right) + \exp\left(\frac{K M^{2} n_{\nu}\xi_{N_{\nu}}^{2}}{2\gamma}\right)\Pi\left(\mathcal{P}\backslash\mathcal{P}_{N_{\nu}}\right)\nonumber\\
& + \exp\left(\frac{K M^{2} n_{\nu}\xi_{N_{\nu}}^{2}}{2\gamma}\right)\cdot\mathbb{E}_{P_{0},P_{\nu}}\displaystyle\mathop{\int}_{\left\{P\in\mathcal{P}_{N_{\nu}}:d^{\pi}_{N_{\nu}}\left(P,P_{0}\right)\geq M\xi_{N_{\nu}}\right\}}\mathop{\prod}_{i=1}^{N_{\nu}}\frac{p^{\pi}}{p_{0}^{\pi}}\left(X_{i}\delta_{\nu i}\right)d\Pi\left(P\right)\left(1-\phi_{n_{\nu}}\right)\nonumber\\
&\mathop{\leq}^{(i)}  \frac{16\gamma^{2}\left[\gamma+C_{3}\right]}{\left(K M^{2}f - 2\gamma\right)^{2} n_{\nu}\xi_{N_{\nu}}^{2}} + \exp\left(-\frac{K M^{2} n_{\nu}\xi_{N_{\nu}}^{2}}{2\gamma}\right) \nonumber\\
&+
\exp\left(\frac{K M^{2} n_{\nu}\xi_{N_{\nu}}^{2}}{2\gamma}\right)\cdot\mathbb{E}_{P_{0},P_{\nu}}\displaystyle\mathop{\int}_{\left\{P\in\mathcal{P}_{N_{\nu}}:d^{\pi}_{N_{\nu}}\left(P,P_{0}\right)\geq M\xi_{N_{\nu}}\right\}}\mathop{\prod}_{i=1}^{N_{\nu}}\frac{p^{\pi}}{p_{0}^{\pi}}\left(X_{i}\delta_{\nu i}\right)d\Pi\left(P\right)\left(1-\phi_{n_{\nu}}\right)\label{postexp},
\end{align}
where in $(i)$ we have used condition~\ref{sizespace} that bounds from above $\Pi\left(\mathcal{P}\backslash\mathcal{P}_{N_{\nu}}\right)$, the prior mass assigned on the portion of the model space that lies outside the sieve, and have plugged in for constant, $C$.

By conditions~\nameref{existtests}, ~\nameref{bounded} and Lemma~\ref{numerator},
\begin{align}
&\mathbb{E}_{P_{0},P_{\nu}}\displaystyle\mathop{\int}_{\left\{P\in\mathcal{P}_{N_{\nu}}:d^{\pi}_{N_{\nu}}\left(P,P_{0}\right)\geq M\xi_{N_{\nu}}\right\}}\mathop{\prod}_{i=1}^{N_{\nu}}\frac{p^{\pi}}{p_{0}^{\pi}}\left(X_{i}\delta_{\nu i}\right)d\Pi\left(P\right)\left(1-\phi_{N_{\nu}}\right)\nonumber\\
&\leq 2\gamma\exp\left(\frac{-K M^{2}n_{\nu}\xi_{N_{\nu}}^{2}}{\gamma}\right)\nonumber
\end{align}

Returning to the expectation in Equation~\ref{postexp},
\begin{align}
&\mathbb{E}_{P_{0},P_{\nu}}\Pi^{\pi}\left(P\in\mathcal{P}:d^{\pi}_{N_{\nu}}\left(P,P_{0}\right)\geq M\xi_{N_{\nu}}\middle\vert\mathbf{X}_{N_{\nu}}\bm{\delta}_{N_{\nu}}\right)\left(1-\phi_{n_{\nu}}\right)\nonumber\\
&\leq \frac{16\gamma^{2}\left[\gamma+C_{3}\right]}{\left(K M^{2} - 2\gamma\right)^{2} N_{\nu}\xi_{N_{\nu}}^{2}} +
\exp\left(-\frac{K M^{2}n_{\nu}\xi^{2}_{N_{\nu}}}{2\gamma}\right) + \exp\left(\frac{K M^{2} n_{\nu}\xi_{N_{\nu}}^{2}}{2\gamma}\right)\cdot2\gamma\exp\left(-\frac{K M^{2}n_{\nu}\xi_{N_{\nu}}^{2}}{\gamma}\right)\nonumber\\
&\mathop{\leq}^{(i)} \frac{16\gamma^{2}\left[\gamma+C_{3}\right]}{\left(K f - 2\gamma\right)^{2}N_{\nu}\xi_{N_{\nu}}^{2}} + 3\gamma\exp\left(-\frac{K M^{2}n_{\nu}\xi_{N_{\nu}}^{2}}{2\gamma}\right),
\end{align}
where in $(i)$ we use our earlier stated bound, $KM^{2}-1> K \rightarrow KM^{2} > K + 1$.

Bringing all the pieces together,
\begin{align}\label{alltogether}
&\mathbb{E}_{P_{0},P_{\nu}}\Pi^{\pi}\left(P\in\mathcal{P}:d^{\pi}_{N_{\nu}}\left(P,P_{0}\right)\geq M\xi_{N_{\nu}}\middle\vert X_{1}\delta_{\nu 1},\ldots,X_{N_{\nu}}\delta_{N_{\nu}}\right)\nonumber\\
&\leq  2\exp\left(-K n_{\nu}\xi_{N_{\nu}}^{2}\right)+  \frac{16\gamma^{2}\left[\gamma+C_{3}\right]}{\left(K f - 2\gamma\right)^{2}N_{\nu}\xi_{N_{\nu}}^{2}} + 3\gamma\exp\left(-\frac{K M^{2}n_{\nu}\xi_{N_{\nu}}^{2}}{2\gamma}\right)\nonumber\\
&\leq\frac{16\gamma^{2}\left[\gamma+C_{3}\right]}{\left(K f - 2\gamma\right)^{2}N_{\nu}\xi_{N_{\nu}}^{2}} + 5\gamma\exp\left(-\frac{K n_{\nu}\xi_{N_{\nu}}^{2}}{2\gamma}\right)
\end{align}
where $\gamma \geq 1$ and $C_{3}> 0$. The right-hand side of Equation~\ref{alltogether} tends to $0$ (as $\nu\uparrow\infty$) in $P_{0}$ probability.
This concludes the proof.
\end{proof}

\section{Enabling Lemmas} \label{AppEnabling}
We next construct two enabling results needed to prove Theorem~\ref{main} to account informative sampling under \nameref{bounded}, \nameref{independence} and \nameref{fraction}.  The first enabling result, Lemma~\ref{numerator}, extends the applicability of \citet{ghosal2007} - Lemmas $2$ and $9$ for $inid$ models to informative sampling without replacement.  This result is used to bound from above the numerator for the expectation with respect to the joint distribution for population generation and the taking of the informative sample, $\left(P_{0},P_{\nu}\right)$, of the pseudo posterior distribution in Equation~\ref{inform_post} on the restricted set of measures, $\{P \in B\}$, where $B = \left\{P\in\mathcal{P}:d_{N_{\nu}}\left(P,P_{0}\right)> \delta\xi_{N_{\nu}}\right\}$, (for any $\delta > 0$).  The restricted set includes those $P$ that are at some \emph{minimum} distance, $\delta\xi_{N_{\nu}}$, from $P_{0}$ under pseudo Hellinger metric, $d^{\pi}_{N_{\nu}}$.  The second result, Lemma~\ref{denominator}, extends Lemma $8.1$ of \citet{Ghosal00convergencerates} to bound the probability of the denominator of Equation~\ref{inform_post} with respect to $\left(P_{0},P_{\nu}\right)$, from below.

\begin{lemma}\label{numerator}
Suppose conditions~\nameref{existtests} and ~\nameref{bounded} hold.  Then for every $\xi > \xi_{N_{\nu}}$, a constant, $K>0$, and any constant, $\delta > 0$,
\begin{eqnarray}
\mathbb{E}_{P_{0},P_{\nu}}\left[\mathop{\int}_{P\in\mathcal{P}\backslash\mathcal{P}_{N_{\nu}}}\mathop{\prod}_{i=1}^{N_{\nu}}
\frac{p^{\pi}}{p_{0}^{\pi}}\left(\mbf{X}_{i}\delta_{\nu i}\right)d\Pi\left(P\right)\left(1-\phi_{n_{\nu}}\right)\right] &\leq& \Pi\left(\mathcal{P}\backslash\mathcal{P}_{N_{\nu}}\right) \label{outside}\\
\mathbb{E}_{P_{0},P_{\nu}}\left[\mathop{\int}_{P\in\mathcal{P}_{N_{\nu}}:d^{\pi}_{N_{\nu}}\left(P,P_{0}\right)> \delta\xi}\mathop{\prod}_{i=1}^{N_{\nu}}
\frac{p^{\pi}}{p_{0}^{\pi}}\left(\mbf{X}_{i}\delta_{\nu i}\right)d\Pi\left(P\right)\left(1-\phi_{n_{\nu}}\right)\right] &\leq& 2\gamma\exp\left(\frac{-K n_{\nu}\delta^{2}\xi^{2}}{\gamma}\right).\label{inside}
\end{eqnarray}
\end{lemma}

The constant multiplier, $\gamma \geq 1$, defined in condition~\nameref{bounded}, restricts the distribution of the sampling design by bounding all marginal inclusion probabilities for population units away from $0$.  As with the main result, the upper bound is injured by $\gamma$.

\begin{proof}\label{AppNumerator}
We proceed constructively to simplify the form of the expectations on the left-hand side of both Equations~\ref{outside} and \ref{inside} and follow with an application of Lemma 2 (and result $2.2$) and Lemma 9 of \citet{ghosal2007}, which is used to establish the right-hand bound of Equation~\ref{inside} (based on the existence of tests, $\phi_{n_{\nu}}$).

Fixing $\nu$, we index units that comprise the population with, $U_{\nu} = \left\{1,\ldots,N_{\nu}\right\}$.  Next, draw a single observed sample of $n_{\nu}$ units from $U_{\nu}$, indexed by subsequence, \newline$\left\{i_{\ell} \in U_{\nu}: \delta_{\nu i_{\ell}} = 1,~\ell = 1,\ldots,n_{\nu}\right\}$.  Without loss of generality, we simplify notation to follow by indexing the observed sample, sequentially, with $\ell = 1,\ldots,n_{\nu}$.

We next decompose the expectation under the joint distribution with respect to population generation, $P_{0}$, and the drawing of a sample, $P_{\nu}$,

Suppose we draw $P$ from some set $B \subset \mathcal{P}$.  By Fubini,
\begin{align}
&\mathbb{E}_{P_{0},P_{\nu}}\left[\mathop{\int}_{P\in B}\mathop{\prod}_{i=1}^{N_{\nu}}
\frac{p^{\pi}}{p_{0}^{\pi}}\left(\mbf{X}_{i}\delta_{\nu i}\right)d\Pi\left(P\right)\left(1-\phi_{n_{\nu}}\right)\right] \nonumber\\
&\leq \mathop{\int}_{P\in B}\left[\mathbb{E}_{P_{0},P_{\nu}}\mathop{\prod}_{i=1}^{N_{\nu}}
\frac{p^{\pi}}{p_{0}^{\pi}}\left(\mbf{X}_{i}\delta_{\nu i}\right)\left(1-\phi_{n_{\nu}}\right)\right]d\Pi\left(P\right)\\
&\leq \mathop{\int}_{P\in B}
\left\{\displaystyle\mathop{\sum}_{\bm{\delta}_{\nu}\in\Delta_{\nu}}\mathbb{E}_{P_{0}}\left[\mathop{\prod}_{\ell=1}^{n_{\nu}}
\left[\frac{p}{p_{0}}\left(\mbf{X}_{\ell}\right)\right]^{\frac{1}{\pi_{\nu \ell}}}\left(1-\phi_{n_{\nu}}\right)\middle\vert \bm{\delta}_{\nu}\right]P_{P_{\nu}}\left(\bm{\delta}_{\nu}\right)\right\}d\Pi\left(P\right)\\
&\leq \mathop{\int}_{P\in B}
\mathop{\max}_{\bm{\delta}_{\nu}\in\Delta_{\nu}}\mathbb{E}_{P_{0}}\left[\mathop{\prod}_{\ell=1}^{n_{\nu}}
\left[\frac{p}{p_{0}}\left(\mbf{X}_{\ell}\right)\right]^{\frac{1}{\pi_{\nu \ell}}}\left(1-\phi_{n_{\nu}}\right)\middle\vert \bm{\delta}_{\nu}\right]d\Pi\left(P\right)\\
&\leq \mathop{\int}_{P\in B}
\mathbb{E}_{P_{0}}\left[\mathop{\prod}_{\ell=1}^{n_{\nu}}
\left[\frac{p}{p_{0}}\left(\mbf{X}_{\ell}\right)\right]^{\frac{1}{\pi_{\nu \ell}}}\left(1-\phi_{n_{\nu}}\right)\middle\vert \bm{\delta}^{\ast}_{\nu}\right]d\Pi\left(P\right)\\
&\leq \mathop{\int}_{P\in B}
\mathbb{E}_{P_{0}}\left[\mathop{\prod}_{\ell=1}^{n_{\nu}}
\left[\frac{p}{p_{0}}\left(\mbf{X}_{\ell}\right)\right]\left(1-\phi_{n_{\nu}}\right)\middle\vert \bm{\delta}^{\ast}_{\nu}\right]d\Pi\left(P\right)\label{ratio}\\
&\leq\mathop{\int}_{P\in B}P_{\bm{\delta}^{\ast}_{\nu}}\left(1-\phi_{n_{\nu}}\right)d\Pi\left(P\right)\nonumber,
\end{align}
where $\displaystyle\mathop{\sum}_{\bm{\delta}_{\nu}\in\Delta_{\nu}}P_{P_{\nu}}\left(\bm{\delta}_{\nu}\right) = 1$ \citep{model2003sarndal} and
$\bm{\delta}^{\ast}_{\nu}\in \Delta_{\nu} = \left\{\left\{\delta^{\ast}_{\nu i}\right\}_{i = 1,\ldots,N_{\nu}},~\delta^{\ast}_{\nu i}\in\{0,1\} \right\}$ denotes that sample, drawn from the space of all possible samples, $\Delta_{\nu}$, which maximizes the probability under the population generating distribution for the event of interest.  The inequality in Equation~\ref{ratio} results from $\frac{p}{p_{0}}\leq 1$ and $\frac{1}{\pi_{\nu\ell}}\geq 1$.  The conditional expectation of $\left(1-\phi_{n_{\nu}}\right)$ given $\bm{\delta}^{\ast}_{\nu}$ is denoted by, $P_{\bm{\delta}^{\ast}_{\nu}}\left(1-\phi_{n_{\nu}}\right)$.

If $P \in \mathcal{P}\backslash\mathcal{P}_{N_{\nu}}$,
\begin{align*}
&\mathbb{E}_{P_{0},P_{\nu}}\left[\mathop{\int}_{P\in \mathcal{P}\backslash\mathcal{P}_{N_{\nu}}}\mathop{\prod}_{i=1}^{N_{\nu}}
\frac{p^{\pi}}{p_{0}^{\pi}}\left(\mbf{X}_{i}\delta_{\nu i}\right)\left(1-\phi_{n_{\nu}}\right)\right]d\Pi\left(P\right)\\
&\leq \mathop{\int}_{P\in \mathcal{P}\backslash\mathcal{P}_{N_{\nu}}}P_{\bm{\delta}^{\ast}_{\nu}}\left(1-\phi_{n_{\nu}}\right)d\Pi\left(P\right)
\leq \mathop{\int}_{P \in\mathcal{P}\backslash\mathcal{P}_{N_{\nu}}}d\Pi\left(P\right) = \Pi\left(\mathcal{P}\backslash\mathcal{P}_{N_{\nu}}\right),
\end{align*}
since $\left(1-\phi_{n_{\nu}}\right) \leq 1$.

We next establish a bound for $P_{\bm{\delta}^{\ast}_{\nu}}\left(1-\phi_{n_{\nu}}\right)$ on a sieve or slice. Let $\mathcal{A}^{\pi}_{r} =\{P\in \mathcal{P}_{N_{\nu}}: r\epsilon_{N_{\nu}} \leq d_{N_{\nu}}^{\pi}\left(P,P_{0}\right) \leq 2r\epsilon_{N_{\nu}} \}$ for integers, $r$.  Under observed $\displaystyle\left(\mbf{X}_{1}\delta^{\ast}_{\nu 1},\ldots,\mbf{X}_{N_{\nu}}\delta^{\ast}_{\nu N_{\nu}}\right) \in \mathcal{X}$, by conditions ~\nameref{existtests} and ~\nameref{bounded} we have,

\begin{align}
&\mathop{\sup}_{P\in\mathcal{A}^{\pi}_{r}}P_{\bm{\delta}^{\ast}_{\nu}}\left(1-\phi_{n_{\nu}}\right)\\
&= \sup_{\{P\in \mathcal{P}_{N_{\nu}}:r\xi \leq d_{N_{\nu}}^{\pi}\left(P,P_{0}\right) \leq 2r\xi\}}P_{\bm{\delta}^{\ast}_{\nu}}\left(1-\phi_{n_{\nu}}\right)\\
&\mathop{\leq}^{(i)}\sup_{\left\{P\in \mathcal{P}_{N_{\nu}}:\frac{r\xi}{\sqrt{\gamma}} \leq d_{N_{\nu}}\left(P,P_{0}\right) \leq \frac{2r\xi}{\sqrt{\gamma}}\right\}}P_{\bm{\delta}^{\ast}_{\nu}}\left(1-\phi_{n_{\nu}}\right)\\
&\mathop{\leq}^{(ii)}\exp\left(-\frac{K n_{\nu} r^{2}\xi^{2}}{\gamma}\right),
\end{align}
where the smaller range in $(i)$, $P\in \mathcal{P}_{N_{\nu}}:\frac{r\xi}{\sqrt{\gamma}} \leq d_{N_{\nu}}\left(P,P_{0}\right) \leq \frac{2r\xi}{\sqrt{\gamma}}$, increases $P_{\bm{\delta}^{\ast}_{\nu}}\left(1-\phi_{n_{\nu}}\right)$. The result in $(ii)$ uses condition~\nameref{sizespace} to obtain the result of Lemmas 2 and 9 in \citet{ghosal2007} where we set $\xi \rightarrow \xi/\sqrt{\gamma}$.

Finally, fixing some value for $\delta>0$, set $r = 2^{\ell}\delta$ for a given, for integers, $\ell \geq 0$.  Following the approach for bounding the sum over the slices in \citet{wong1995}, let $L$ be the smallest integer such that $2^{2L}\delta^{2}\xi^{2} > 2\gamma$, since $d_{N_{\nu}}^{\pi} < \sqrt{2\gamma}$ (by our definition of the pseudo Hellinger metric in Section~\ref{empirical}).  Then,
\begin{align}
&\mathbb{E}_{P_{\theta_{0}},P_{\nu}}\left[\mathop{\int}_{\left\{P\in \mathcal{P}_{N_{\nu}}:d_{N_{\nu}}^{\pi}\left(P,P_{0}\right) \geq \delta\xi\right\}}\displaystyle\mathop{\prod}_{i=1}^{N_{\nu}}\frac{p^{\pi}}{p_{0}^{\pi}}\left(\mathbf{X}_{i}\delta_{\nu i}\right)d\Pi\left(P\right)\left(1-\phi_{n_{\nu}}\right)\right]\\
&=\mathop{\sum}_{\ell = 0}^{L}\mathbb{E}_{P_{\theta_{0}},P_{\nu}}\int_{\left\{P\in \mathcal{P}_{N_{\nu}}:2^{\ell}\delta\xi \leq d_{N_{\nu}}^{\pi}\left(P,P_{0}\right) \leq 2^{\ell+1}\delta\xi\right\}}\displaystyle\mathop{\prod}_{i=1}^{N_{\nu}}\frac{p^{\pi}}{p_{0}^{\pi}}\left(\mathbf{X}_{i}\delta_{\nu i}\right)d\Pi\left(P\right)\left(1-\phi_{N_{\nu}}\right)\\
&\leq\gamma\mathop{\sum}_{\ell = 0}^{L}\exp\left(-\frac{2^{2\ell}K n_{\nu}\delta^{2}\xi^{2}}{\gamma}\right)\\
&\leq 2\gamma\exp\left(-\frac{K n_{\nu}\delta^{2}\xi^{2}}{\gamma}\right),
\end{align}
for $n_{\nu}$ sufficiently large such that $\frac{K n_{\nu}\delta^{2}\xi^{2}}{\gamma}\geq 1$.

This concludes the proof.

\end{proof}

\begin{lemma}\label{denominator}
For every $\xi > 0$ and measure $\Pi$ on the set,
\begin{equation*}
B = \left\{P:-P_{0}\log\left(\frac{p}{p_{0}}\right) \leq \xi^2, P_{0}\left(\log\frac{p}{p_{0}}\right)^{2} \leq \xi^{2}\right\}
\end{equation*}
under the conditions ~\nameref{sizespace}, ~\nameref{priortruth}, ~\nameref{bounded}, and ~\nameref{independence}, we have for every $C > 0 $ and $N_{\nu}$ sufficiently large,
\begin{equation}\label{denomresult}
\mbox{Pr}\left\{\mathop{\int}_{P\in\mathcal{P}}\displaystyle\mathop{\prod}_{i=1}^{N_{\nu}}\frac{p^{\pi}}{p_{0}^{\pi}}
\left(\mbf{X}_{i}\delta_{\nu i}\right)d\Pi\left(P\right)\leq \exp\left[-(1+C)N_{\nu}\xi^{2}\right]\right\}
\leq \frac{\gamma+C_{3}}{C^{2} N_{\nu}\xi^{2}},
\end{equation}
where the above probability is taken with the respect to $P_{0}$ and the sampling generating distribution, $P_{\nu}$, jointly.
\end{lemma}

The bound of ``$1$" in the numerator of the result for Lemma $8.1$ of \citet{Ghosal00convergencerates}, is replaced with $\gamma + C_{3}$ for our generalization of this result in Equation~\ref{denomresult}.  The sum of positive constants, $\gamma + C_{3}$, is greater than $1$ and will be larger for sampling designs where the inclusion probabilities, $\{\pi_{\nu i}\}$, express relatively higher gradients.  Observing each finite population in a skewed fashion through the taking of an informative sample may only slow the rate of posterior contraction (as compared to contraction of the posterior distribution defined on the \emph{fully} observed finite population).

\begin{proof}\label{AppDenominator}
By Jensen's inequality,
\begin{align*}
\log\mathop{\int}_{P\in\mathcal{P}}\mathop{\prod}_{i=1}^{N_{\nu}}\frac{p^{\pi}}{p_{0}^{\pi}}\left(\mbf{X}_{i}\delta_{\nu i}\right)d\Pi\left(P\right) &\geq \mathop{\sum}_{i=1}^{N_{\nu}}\displaystyle\mathop{\int}_{P\in\mathcal{P}}\log\frac{p^{\pi}}{p_{0}^{\pi}}\left(\mbf{X}_{i}\delta_{\nu i}\right)d\Pi\left(P\right)\\
&= N_{\nu}\cdot\mathbb{P}_{N_{\nu}}\mathop{\int}_{P\in\mathcal{P}}\log\frac{p^{\pi}}{p_{0}^{\pi}}d\Pi\left(P\right),
\end{align*}
where we recall that the last equation denotes the empirical expectation functional taken with respect to the joint distribution over population generating and informative sampling.  By Fubini,
\begin{align*}
\mathbb{P}_{N_{\nu}}\mathop{\int}_{P\in\mathcal{P}}\log\frac{p^{\pi}}{p_{0}^{\pi}}d\Pi\left(P\right)
&= \mathop{\int}_{P\in\mathcal{P}}\left[\mathbb{P}_{N_{\nu}}\log\frac{p^{\pi}}{p_{0}^{\pi}}\right]d\Pi\left(P\right)\\
&= \mathop{\int}_{P\in\mathcal{P}}\left[\mathbb{P}_{N_{\nu}}\frac{\delta_{\nu}}{\pi_{\nu}}
\log\frac{p}{p_{0}}\right]d\Pi\left(P\right)\\
&= \mathop{\int}_{P\in\mathcal{P}}\left[\mathbb{P}^{\pi}_{N_{\nu}}
\log\frac{p}{p_{0}}\right]d\Pi\left(P\right)\\
&= \mathbb{P}^{\pi}_{N_{\nu}}\mathop{\int}_{P\in\mathcal{P}}\log\frac{p}{p_{0}}d\Pi\left(P\right),
\end{align*}
where we, again, apply Fubini.

Then, the probability statement in the result of Equation~\ref{denomresult} is bounded (from above) by,
\begin{align*}
&\mbox{Pr}\left\{\mathbb{G}^{\pi}_{N_{\nu}}\mathop{\int}_{P\in\mathcal{P}}\log\frac{p}{p_{0}}
d\Pi\left(P\right) \leq -\sqrt{N_{\nu}}\xi^{2}\left(1+C\right)
-\sqrt{N_{\nu}}P_{0}\mathop{\int}_{P\in\mathcal{P}}\log\frac{p}{p_{0}}d\Pi\left(P\right)\right\}\\
&=\mbox{Pr}\left\{\mathbb{G}^{\pi}_{N_{\nu}}\mathop{\int}_{P\in\mathcal{P}}\log\frac{p}{p_{0}}
d\Pi\left(P\right)\leq -\sqrt{N_{\nu}}\xi^{2}\left(1+C\right)
-\sqrt{N_{\nu}}\mathop{\int}_{P\in\mathcal{P}}P_{0}\log\frac{p}{p_{0}}d\Pi\left(P\right)\right\}\\
&=\mbox{Pr}\left\{\mathbb{G}^{\pi}_{N_{\nu}}\mathop{\int}_{P\in\mathcal{P}}\log\frac{p}{p_{0}}
d\Pi\left(P\right)\leq -\sqrt{N_{\nu}}\xi^{2}\left(1+C\right) + \sqrt{N_{\nu}}\xi^{2} = -\sqrt{N_{\nu}}\xi^{2}C\right\},
\end{align*}
where we have again applied Fubini in the second inequality and also the bound for $P_{0}\log\frac{p}{p_{0}}\leq\xi^{2}$
for $P$ on the set $B$.

We now apply Chebyshev and Jensen's inequality to bound the probability,
\begin{subequations}\label{chebyshev}
\begin{align}
\mbox{Pr}\left\{\mathbb{G}^{\pi}_{N_{\nu}}\mathop{\int}_{P\in\mathcal{P}}\log\frac{p}{p_{0}}
d\Pi\left(P\right)\leq -\sqrt{N_{\nu}}\xi^{2}C\right\}
&\leq\frac{\Var\left[\mathop{\int}_{P\in\mathcal{P}}\mathbb{G}^{\pi}_{N_{\nu}}\log\frac{p}{p_{0}}d\Pi\left(P\right)\right]}
{N_{\nu}\xi^{4}C^{2}}\\
&\leq\frac{\displaystyle\mathop{\int}_{P\in\mathcal{P}}\left[\Var\left(\mathbb{G}^{\pi}_{N_{\nu}}\log\frac{p}{p_{0}}\right)\right]d\Pi\left(P\right)}
{N_{\nu}\xi^{4}C^{2}}\label{chebyshev:var}\\
&\leq\frac{\displaystyle\mathop{\int}_{P\in\mathcal{P}}\left[\mathbb{E}_{P_{0},P_{\nu}}
\left(\mathbb{G}^{\pi}_{N_{\nu}}\log\frac{p}{p_{0}}\right)^{2}\right]d\Pi\left(P\right)}
{N_{\nu}\xi^{4}C^{2}}\label{chebyshev:e2}\\
&\leq\frac{\displaystyle\mathop{\int}_{P\in\mathcal{P}}\left[\mathbb{E}_{P_{0},P_{\nu}}
\left(\sqrt{N_{\nu}}\mathbb{P}^{\pi}_{N_{\nu}}\log\frac{p}{p_{0}}\right)^{2}\right]d\Pi\left(P\right)}
{N_{\nu}\xi^{4}C^{2}}\label{chebyshev:p},
\end{align}
\end{subequations}
where $\mathbb{E}_{P_{0},P_{\nu}}\left(\cdot\right)$ denotes the expectation with respect to the joint distribution over population generation and sampling (from that population) without replacement.   We apply Jensen's inequality in Equation~\ref{chebyshev:var} and use $\mathbb{E}\left(X^{2}\right) > \Var\left(X\right)$ in the third inequality, stated in Equation~\ref{chebyshev:e2}, and drop the centering term in Equation~\ref{chebyshev:p}.  We now bound the expectation inside the square brackets on the right-hand side of Equation~\ref{chebyshev:p}, which is taken with respect to this joint distribution.  In the sequel, define $\mathcal{A}_{\nu} = \sigma\left(\mbf{X}_{1},\ldots,\mbf{X}_{N_{\nu}}\right)$ as the sigma field of information potentially available for the $N_{\nu}$ units in population, $U_{\nu}$.

\begin{align*}
\mathbb{E}_{P_{0},P_{\nu}}
\left(\sqrt{N_{\nu}}\mathbb{P}^{\pi}_{N_{\nu}}\log\frac{p}{p_{0}}\right)^{2} &=
\frac{1}{N_{\nu}}\mathop{\sum}_{i,j\in U_{\nu}}\mathbb{E}_{P_{0},P_{\nu}}\left(\frac{\delta_{\nu i}\delta_{\nu j}}
{\pi_{\nu i}\pi_{\nu j}}\log\frac{p}{p_{0}}\left(\mbf{X}_{i}\right)\log\frac{p}{p_{0}}\left(\mbf{X}_{j}\right)\right)\\
&= \displaystyle\frac{1}{N_{\nu}}\mathop{\sum}_{i=j\in U_{\nu}}\mathbb{E}_{P_{0}}\left[\mathbb{E}_{P_{\nu}}\left\{\left(\frac{\delta_{\nu i}}
{\pi_{\nu i}^{2}}\left(\log\frac{p}{p_{0}}\left(\mbf{X}_{i}\right)\right)^{2}\right)\middle\vert \mathcal{A}_{\nu}\right\}\right]\\
&+ \frac{1}{N_{\nu}^{2}}\mathop{\sum}_{i\neq j\in U_{\nu}}\mathbb{E}_{P_{0}}\left[\frac{\mathbb{E}_{P_{\nu}}\left[\delta_{\nu i}\delta_{\nu j}\vert\mathcal{A}_{\nu}
\right]}{\pi_{\nu i}\pi_{\nu j}}\log\frac{p}{p_{0}}\left(\mbf{X}_{i}\right)\log\frac{p}{p_{0}}\left(\mbf{X}_{j}\right)\right]\\
&= \displaystyle\frac{1}{N_{\nu}}\mathop{\sum}_{i=j\in U_{\nu}}\mathbb{E}_{P_{0}}\left[\left(\frac{1}
{\pi_{\nu i}}\right)\left(\log\frac{p}{p_{0}}\left(\mbf{X}_{i}\right)\right)^{2}\right]\\
&+ \frac{1}{N_{\nu}}\mathop{\sum}_{i\neq j\in U_{\nu}}\mathbb{E}_{P_{0}}\left[\frac{\pi_{\nu ij}}{\pi_{\nu i}\pi_{\nu j}}\log\frac{p}{p_{0}}\left(\mbf{X}_{i}\right)\log\frac{p}{p_{0}}\left(\mbf{X}_{j}\right)\right]\\
&\displaystyle\leq\xi^{2}\mathop{\sup}_{\nu}\left[\frac{1}{\displaystyle\mathop{\min}_{i\in U_{\nu}}\pi_{\nu i}}\right] + \xi^{2}\left(N_{\nu}-1\right)\mathop{\sup}_{\nu}\displaystyle\mathop{\max}_{i \neq j\in U_{\nu}}\left[\displaystyle\middle\vert\frac{\pi_{\nu ij}}{\pi_{\nu i}\pi_{\nu j}}\middle\vert\right] \\
&\leq\xi^{2}\left(\gamma + C_{3}\right),
\end{align*}
for sufficiently large $N_{\nu}$, where we have applied the condition for $P\in B$ for the first term of the last two inequalities and conditions and ~\nameref{bounded} and ~\nameref{independence} for the last inequality.  We additionally note that $\pi_{\nu ij} = \pi_{\nu j}$ when $i = j,~i,j \in U_{\nu}$. This concludes the proof.
\end{proof}

\end{document}